\documentclass[11pt]{article}

\usepackage{amsmath}
\usepackage{mathtools}
\usepackage{amssymb}
\usepackage{latexsym}
\usepackage{enumitem}
\usepackage{mathrsfs}
\usepackage{comment}
\usepackage{pifont}
\usepackage{enumitem}
\usepackage{bbold}
\usepackage{graphicx}
\usepackage{color}
\usepackage{hyperref}
\hypersetup{
    colorlinks=true,
    linkcolor=blue,
    filecolor=magenta,      
    urlcolor=cyan,
    citecolor=magenta
}

\newtheorem{assumption}{Assumption}

\def\qed{ \ \vrule width.2cm height.2cm depth0cm\smallskip}
\newenvironment{proof}{\noindent {\bf Proof.\/}}{$\qed$\vskip 0.1in}

\newcommand{\la}{\langle}
\newcommand{\ra}{\rangle}

\newcommand{\ol}{\overline}
\newcommand{\ul}{\underline}

\newcommand{\ba}{\begin{array}}
\newcommand{\ea}{\end{array}}
\newcommand{\be}{\begin{equation}}
\newcommand{\ee}{\end{equation}}
\newcommand{\bea}{\begin{eqnarray}}
\newcommand{\eea}{\end{eqnarray}}
\newcommand{\beaa}{\begin{eqnarray*}}
\newcommand{\eeaa}{\end{eqnarray*}}

\newcommand{\essinf}{\operatornamewithlimits{essinf}}

\def\dbD{\mathbb{D}}
\def\dbE{\mathbb{E}}
\def\dbF{\mathbb{F}}

\def\dbI{\mathbb{I}}

\def\dbL{\mathbb{L}}

\def\dbP{\mathbb{P}}
\def\dbR{\mathbb{R}}
\def\dbS{\mathbb{S}}

\def\g{\gamma}
\def\d{\delta}
\def\e{\varepsilon}

\def\k{\kappa}
\def\l{\lambda}
\def\m{\mu}
\def\n{\nu}
\def\si{\sigma}
\def\t{\tau}
\def\f{\varphi}
\def\th{\theta}
\def\o{\omega}

\def\D{\Delta}

\def\L{\Lambda}

\def\O{\Omega}
\def\cA{{\cal A}}
\def\cB{{\cal B}}

\def\cF{{\cal F}}

\def\cL{{\cal L}}

\def\cN{{\cal N}}
\def\cO{{\cal O}}
\def\cP{{\cal P}}

\def\cS{{\cal S}}
\def\cT{{\cal T}}

\def\cW{{\cal W}}

\def\cZ{{\cal Z}}

\def\ch{\textsc{h}}

\def\no{\noindent}

\def\ms{\medskip}

\def\q{\quad}
\def\qq{\qquad}

\def\pa{\partial}
\def\cd{\cdot}
\def\cds{\cdots}

\def\tr{\hbox{\rm tr}}

\def\qed{ \hfill \vrule width.25cm height.25cm depth0cm\smallskip}

\newcommand{\basa}{\begin{assumption}}
\newcommand{\easa}{\end{assumption}}

\newcommand{\bas}{\begin{assum}}
\newcommand{\eas}{\end{assum}}

\def\essinf{\mathop{\rm ess\;inf}}

\def\pa{\partial}

 \def\cd{\cdot}
\def\cds{\cdots}

\def\tr{\hbox{\rm tr$\,$}}

\def\dis{\displaystyle}

\def\bff{{\bf f}}
\def\bx{{\bf x}}
\def\bi{{\bf i}}
\def\be{{\bf e}}
\def\bI{{\bf I}}
\def\bX{{\bf X}}
\def\btau{{\boldsymbol \tau}}

\def\1{{\bf 1}}
\def\by{{\bf y}}
\def\bY{{\bf Y}}

\def\:{\!:\!}
\def\reff#1{{\rm(\ref{#1})}}
\def \proof{{\noindent \bf Proof\quad}}

\def \bS{{\bf S}}

\def \bQ{{\bf Q}}

\DeclareMathAlphabet\mathbfcal{OMS}{cmsy}{b}{n}

at 9pt

\begin{document}

\newtheorem{thm}{Theorem}[section]
\newtheorem{lem}[thm]{Lemma}
\newtheorem{cor}[thm]{Corollary}
\newtheorem{prop}[thm]{Proposition}
\newtheorem{rem}[thm]{Remark}
\newtheorem{eg}[thm]{Example}
\newtheorem{defn}[thm]{Definition}
\newtheorem{assum}[thm]{Assumption}

\numberwithin{equation}{section}

\title{\bf{From finite population optimal stopping 
             \\ to mean field optimal stopping}\footnote{The first two authors are grateful for the financial support from the Chaires FiME-FDD and Financial Risks of the Louis Bachelier Institute. The third author is supported in part by NSF grants DMS-1908665 and DMS-2205972.}}
\author{Mehdi Talbi\footnote{Laboratoire de Probabilités, Statistiques et Modélisation, Université Paris-Cité, France, talbi@lpsm.paris} \quad Nizar Touzi\footnote{Tandon School of Engineering, New York University, United States, nt2635@nyu.edu}  \quad Jianfeng Zhang\footnote{Department of Mathematics, University of Southern California, United States, jianfenz@usc.edu. }}
\date{\today}

\maketitle

\begin{abstract}
This paper analyzes the convergence of the finite population optimal stopping problem towards the corresponding mean field limit. Building on the viscosity solution characterization of the mean field optimal stopping problem of our previous papers \cite{TTZ, TTZ2}, we prove the convergence of the value functions by adapting the Barles-Souganidis \cite{BarSou} monotone scheme method to our context. We next characterize the optimal stopping policies of the mean field problem by the accumulation points of the finite population optimal stopping strategies. In particular, if the limiting problem has a unique optimal stopping policy, then the finite population optimal stopping strategies do converge towards this solution.
As a by-product of our analysis, we provide an extension of the standard propagation of chaos to the context of stopped McKean-Vlasov diffusions.
\end{abstract}

\no{\bf MSC2020.} 60G40, 35Q89, 49N80, 49L25, 65K15.

\vspace{3mm}
\no{\bf Keywords.}  Mean field optimal stopping, obstacle problems, viscosity solutions, propagation of chaos.


\section{Introduction}
This paper is a continuation of our previous papers \cite{TTZ, TTZ2},  where we established the wellposedness, especially in the viscosity sense, for an obstacle equation on Wasserstein space derived from a mean field optimal stopping problem. In this paper we focus on the convergence of the corresponding $N$-player optimal stopping problem, the latter was investigated by Kobylanski, Quenez \& Rouy-Mironescu \cite{KQR} in a more general setting. We shall obtain both the convergence of the value functions and the propagation of chaos in terms of the (approximate) optimal trajectories.

In the context of mean field games, initiated independently by Lasry \& Lions \cite{LL} and Caines, Huang \& Malhamé \cite{CHM}, the convergence of the $N$-player game to the mean field game is one of the core issues in the field. When the mean field equilibrium is unique, typically under certain monotonicity conditions, one may use the popular master equation approach to obtain the convergence, see Cardaliaguet, Delarue, Lasry \& Lions \cite{CDLL}, followed by Bayraktar \& Cohen \cite{BC}, Bayraktar, Cecchin, Cohen \& Delarue \cite{BCCD2}, Cardaliaguet \cite{Cardaliaguet}, Cecchin \& Pelino \cite{CecPel}, Delarue-Lacker-Ramanan \cite{DLR1,DLR2}, Gangbo \& Meszaros \cite{GM}, and Mou \& Zhang \cite{MZ}, to mention a few. When there are multiple equilibria, there are numerous publications on the convergence issue, mainly on the propagation of chaos, see e.g.  Camona \& Delarue \cite{CD0},  Cecchin, Dai Pra, Fisher \& Pelino \cite{CDFP}, Djete \cite{Djete, Djete1}, Doncel, Gast \& Gaujal \cite{DGG}, Feleqi \cite{Feleqi}, Fischer \cite{Fischer}, Fischer \& Silva \cite{FischerSilva}, Lacker \cite{Lacker1, Lacker2}, Lacker \& Flem \cite{LackerF}, Lasry \& Lions \cite{LL}, Lauriere \& Tangpi \cite{LauTan}, Nutz, San Martin \& Tan \cite{NST},  and Possamai \& Tangpi \cite{PT}. We would also like to mention the set valued approach in Iseri \& Zhang \cite{IZ}, where the object is the set of values over all equilibria. 

In the context of mean field controls, an important stream of works focuses on the convergence issue, see e.g.\ Bayraktar \& Chakraborty \cite{BayCha}, Cardaliaguet, Daudin, Jackson \& Souganidis \cite{CDJS}, Cardaliaguet \& Souganidis \cite{CarSou},  Cavagnari, Lisini, Orrieri \& Savare \cite{CLOS}, Cecchin \cite{Cecchin}, Daudin, Delarue \& Jackson \cite{daudin2023optimal}, Djete, Possamai \& Tan  \cite{DPT},  Fischer \& Livieri \cite{FisLiv},  Fornasier, Lisini, Orrieri \& Savare \cite{FLOS}, Germain, Pham \& Warin \cite{GPW}, and Lacker \cite{LackerPropagation}. In particular, by utilizing some strong regularity of the value function, \cite{GPW, CDJS} obtained certain rate of convergence. \cite{daudin2023optimal} also obtained rates of convergence by mollifying the value function of the mean field problem in order to obtain ``almost" classical subsolutions, which are then projected on $\dbR^N$ to derive estimates for the $N$-particle problem. The work closest to ours is Gangbo, Mayorga \& Swiech \cite{GMS}, which uses the uniqueness of viscosity solution (in certain sense) to prove the convergence of the value function in a setting without idiosyncratic but with common noise. We also mention the paper Talbi \cite{Talbi} by the first author which applies the same approach of this paper to finite dimensional approximation of mean field control problems, including the path dependent case.

The first goal of this paper is to apply the Barles \& Sounganids \cite{BarSou} approach to prove the convergence of the value functions. To our best knowledge, this is the first work in the mean field literature to use the Barles-Sounganids approach. Roughly speaking, denoting by $V$ the value function of the mean field optimal stopping problem and $V^N$ the value function of the corresponding  $N$-player problem, we shall prove that
$\underline V := \underset{N\to\infty}{\lim \inf} \ V^N$ and $\overline V := \underset{N\to\infty}{\lim \sup} \ V^N$
are viscosity supersolution and subsolution, respectively, of our mean field obstacle equation. Then, it follows from the comparison principle of viscosity solutions, which was established in our previous paper \cite{TTZ2}, that $\underline V =\overline V$ and hence $ \lim_{N\to\infty} V^N=V$. Note that our viscosity solution approach allows us to deal with value functions which are merely continuous. Consequently, besides the obvious difference that we consider mean field optimal stopping problem instead of mean field control problem, we require weaker regularity conditions for the coefficients than  \cite{CDJS, GPW}.

Our convergence relies on the viscosity solution property of $V^N$ to a finite-dimensional PDE, as deduced from the general solution of \cite{KQR} expressed in terms of a recursive sequence of Snell envelopes. We shall refer to this equation as a cascade obstacle problem due to its remarkable structure.  One subtle issue is the choice of  test functions for viscosity solutions of this finite-dimensional PDE. To adapt to our notion of viscosity solution for our mean field obstacle equation in \cite{TTZ2}, we only require test functions to be tangent to the super/subsolution through the mean, whereas the tangency is pointwise in the standard literature. Our definition is inspired by the the viscosity theory developed for path-dependent PDEs, see e.g. the series of papers by Ekren, Keller, Ren, Touzi \& Zhang \cite{Ekren, EKTZ, ETZ1, ETZ2,  RTZ}.

Our second main result is the propagation of chaos for stopped McKean-Vlasov diffusion (see e.g.\ Sznitman \cite{Snitzman} for the case of classical diffusions). By using the convergence of the value functions, we establish three results: roughly speaking,

$\bullet$ Any optimal strategy of the mean field problem can be approximated by (approximate) optimal strategies of the $N$-player problems;

$\bullet$ Any limit point of  the (approximate) optimal strategies of the $N$-player problem is  optimal for the mean field problem;

 $\bullet$ If the optimal strategy of the mean field problem is unique, then any (approximate) optimal strategies of the $N$-player problem converges to that unique optimal strategy of the mean field problem.

\no These results are consistent with the convergence results for mean field equilibria in the mean field game literature, see e.g.  \cite{Lacker1}. We shall remark that, for mean field games, the $N$-player problems are quite different between closed loop controls and open loop controls. However, for our mean field optimal stopping problem (and for mean field control problems),  the open loop and the closed loop stopping strategies typically lead to the same value function. In particular, for the convenience of establishing the regularity of the value functions, we consider open loop stopping times for the $N$-player problems.

The paper is organized as follows. In Sections \ref{sect-MF} and \ref{sect-N} we introduce the mean field optimal stopping problem and the finite population optimal stopping problem, respectively. In Section \ref{sect-main} we present the main results of the paper. Sections \ref{sect-convergence} and \ref{sec-proof-propagation} are devoted to the proofs of the main results, one on the convergence of the value functions, and the other on the propagation of chaos. 
In Section \ref{sect-Ndimpb} we study the multiple optimal stopping problem in details. 
Finally, we provide some technical proofs in Appendix.

\vspace{3mm}
\no {\bf Notations.} We denote by $\cP(\O,\cF)$ the set of probability measures on a measurable  space $(\O,\cF)$, and $\cP_2(\O,\cF)$ the subset of probability measures in $\cP(\O,\cF)$ with finite second order moment, equipped with the $2$-Wasserstein distance $\cW_2$. When $\O$ is a topological space and $\cF$ is the Borel $\si$-field, we simply denote them as $\cP(\O)$ and $\cP_2(\O)$. For a random variable $Z$ on a probability space $(\O, \cF, \dbP)$, we denote by $\dbP_Z:=\dbP\circ Z^{-1}$ the law of $Z$ under $\dbP$. The space of $d\times d-$symmetric matrice is denoted by $\cS_d$, with $\cS_d^+$ the subset of non-negative matrices. For vectors $x, y\in \dbR^d$ and matrices $A, B\in \dbR^{d\times N}$, denote $ x\cd y:=\sum_{i=1}^n x_iy_i$  and $A:B:= \tr(A B^\top)$. We shall also write “USC" (resp. “LSC") “upper (resp. lower) semi-continuous". Moreover, denote
 \beaa
\left.\ba{c}
\dis \bS:= \dbR^d \times \{0,1\},\q  \mathbf{Q}_t:=[t,T)\times\cP_2(\mathbf{S}), \q
 \overline{\mathbf{Q}}_t:=[t,T]\times\cP_2(\mathbf{S}),\q t\in [0, T);\ms\\
 \dis  \|m\|_p^p := \int_\bS |x|^p m(dy),~ m\in \cP_2(\bS),\q\mbox{for any}~p\ge 1;
\ea\right. 
 \eeaa
and, for any $N\ge 1$,
\bea
\label{QN}
 \left.\ba{c}
\dis \L^N_t:=[t,T)\times \bS^N, \q \overline{\L}^N_t:=[t,T]\times\bS^N;\ms\\
 \dis m^N(\by) := \frac 1 N \sum_{k=1}^N \d_{y_k} \in \cP(\dbR^d),\q\mbox{for}\q \by =(y_1,\cds, y_N) \in \bS^N;\ms\\
\dis \cP^N(\bS):= \{m_N(\by):  \by\in \bS^N\},\q  \mathbf{Q}^N_t:= [t, T)\times \cP^N(\bS), \q
  \ol \bQ^N_t:=[t, T]\times \cP^N(\bS);\ms\\
\dis \|\bx\|_p^p := {1\over N} \sum_{k=1}^N |x_k|^p,~ \bx=(x_1,\cds, x_N)\in \dbR^{d\times N},
 \ea\right.
 \eea

\section{The mean field optimal stopping problem}
\label{sect-MF}
In this section we introduce the mean field optimal stopping problem and the related obstacle equation investigated in our previous works \cite{TTZ, TTZ2}.
\subsection{The mean field optimal stopping problem}

Let $T < \infty$, and denote $\O := C^0([-1,T],\dbR^d) \times \dbI^0([-1,T])$ the canonical space, with: 

$\bullet$ $C^0([-1,T],\dbR^d)$ is the set of continuous paths from $[-1,T]$ to $\dbR^d$, constant on $[-1,0)$; 

$\bullet$ $\dbI^0([-1,T])$ is the set of non-increasing and càdlàg maps from $[-1,T]$ to $\{0,1\}$,  constant 

~ on $[-1,0)$, and ending with value $0$ at $T$.
 
  \no We recall that $\O$ is Polish under the Skorokhod distance.
 We also denote $Y:= (X,I)$ the canonical process taking values in the state space $\mathbf{S}$, $\dbF = (\cF_t)_{t \in [-1,T]}$ its canonical filtration, and the corresponding jump time of the survival process $I$:
\bea
\label{tau}
\iota := \inf\{t \ge 0 : I_t = 0\},  \ \mbox{so that $I_t := I_{0-}\1_{t < \iota}$ for all $t \in [-1,T]$.}
\eea 
Since $I$ is càdlàg, $\iota$ is an $\dbF-$stopping time. We emphasise on the fact that the extension to $-1$ is arbitrary, as the extension of time to the left of the origin is only needed to allow for an immediate stop at time $t=0$. 

Let $(b,\si,f): [0,T] \times \dbR^d \times \cP_2(\mathbf{S}) \rightarrow \dbR^d \times \dbS_d^+\times \dbR$, and $ g : \cP_2(\dbR^d) \rightarrow \dbR$.  
\begin{assum}
\label{assum-bsig}
\no{\rm (i)} $b, \si$ are continuous in $t$, and uniformly Lipschitz continuous in $(x, m)$.

\no{\rm (ii)} $f$ is Borel measurable and has quadratic growth in $x\in \dbR^d$, and 
\bea
\label{F}
F(t,m) := \int_{\dbR^d} f(t,x,m)m(dx,1)
~\mbox{is continuous on}~[0, T]\times \cP_2(\bS).
\eea
\no{\rm (iii)}  $g$ is upper-semicontinuous and locally bounded.
\end{assum}
 We also extend $g$ to $\cP_2(\bS)$ by $g(m) := g(m(\cd, \{0,1\}))$. Let $V$ be the dynamic value function:
 \begin{equation}\label{weakoptstop}
\dis V(t,m) := \underset{\dbP \in \cP(t,m)}{\sup} J(t,\dbP), ~ (t,m) \in \ol \bQ_0,\q~J(t,\dbP) := \int_t^T  F(r,  \dbP_{Y_r})dr  + g(\dbP_{Y_T}),
 \end{equation}
with $\cP(t,m)$ the set of probability measures $\dbP$ on $(\O,\cF_T)$ s.t.\ $\dbP_{Y_{t-}} = m$ and the processes:
 \bea\label{martingalepb}
M_. := X_. - \int_t^. b(r,X_r,\dbP_{Y_r})I_rdr \q \mbox{and} \q M_. M_.^\intercal -  \int_t^. \si^2(r,X_r,\dbP_{Y_r})I_rdr
\eea
are $\dbP-$martingales on $[t,T]$, that is, for some $\dbP-$Brownian motion $W^\dbP$, 
\begin{equation}\label{asympt}
X_s = X_t + \int_t^s b(r, X_r, \dbP_{Y_r})I_r dr + \sigma(r, X_r, \dbP_{Y_r})I_r dW_r^\dbP , \ I_s = I_{t-} \1_{s < \iota}, \ \dbP-\mbox{a.s.}
\end{equation}
A special element of $\cP(t, m)$ is $\bar \dbP=\bar \dbP^{t,m}$ under which $X$ is unstopped. That is,
\bea
\label{barP}
X_s = X_t + \int_t^s b(r, X_r, \bar \dbP_{Y_r}) I_r dr + \int_t^s \si(r, X_r, \bar \dbP_{Y_r})I_r  dW^{\bar \dbP}_r,~ I_s = I_{t-}\1_{[t, T)}(s),~ \bar \dbP \mbox{-a.s.} 
\eea
Note that $Y_. = Y_{. \wedge \iota}$, and in particular $Y_T = Y_\iota,$ $\dbP-$a.s. Moreover, from the definition of $F$ in \eqref{F}, we have $\int_t^T F(r,\dbP_{Y_r})dr=\dbE^\dbP\int_t^\iota f(r,X_r,\dbP_{Y_r})dr$.

We observe that \eqref{weakoptstop} admits an optimal control as $\cP(t,m)$ is $\cW_2$-compact by \cite{TTZ}.

\subsection{Obstacle equation on Wasserstein space}
We first recall some differential calculus tools on the Wasserstein space. We say that a function $U:\cP_2(\bS) \to \dbR$ has a functional linear derivative $\delta_m U:\cP_2(\bS)\times\bS \to \dbR$ if
 \beaa
U(\tilde m)-U(m) 
= 
\int_0^1 \int_{\bS} \d_m U(l \tilde m + (1-l)m, y)(\tilde m-m)(dy)dl,
\q\mbox{for all}~
m, \tilde m\in\cP_2(\bS),
\eeaa
and $\d_m U(m,\cdot)$ has quadratic growth, locally uniformly in $m$, so as to guarantee integrability in the last expression. As in \cite{TTZ}, we denote
\beaa
\d_m U_i(t,m,x) := \d_m U(t,m,x,i)  &\mbox{for}& i \in \{0,1\},\qq
D_I U := \d_m U_1 - \d_m U_0,
\eeaa
and we introduce the measure flow generator of $X$
\bea\label{DiffOp}
\dbL U(t,m)
&:=& 
\pa_t U(t,m) + \int_{\dbR^d}  \cL_x\d_m U_1(t,m,x) m(dx,1), \nonumber
\eea
where $\cL_x$ is the generator of $X$: 
\bea
\label{cLx}
\cL_x\d_m U_1
&:=&
b\cd \pa_x \d_m U_1 + {1\over 2} \si^2\!\!:\! \pa_{xx}^2  \d_m U_1. \nonumber
\eea
 We may write $\pa_{\m} U_{1,1}(t,m,x) :=  \pa_x \d_m U(t,m,x,1)$ in the spirit of the Lions derivative, and we also introduce $\pa_{\m \m}^2 U_{1,1}(t,m,x,\tilde x) :=  \pa_{\tilde x} \d_m \pa_x \d_m U(t,m,x,1,\tilde x,1)\in \dbR^{d\times d}.$

We next introduce a partial order $\preceq$ on $\cP_2(\mathbf{S})$. We say that $m' \preceq m$ if $m'(dx,1)$ is absolutely continuous with respect to $m(dx,1)$ with density $p : \dbR^d \rightarrow [0,1]$, i.e.,
\bea\label{order}
m'(dx, 1) = p(x) m(dx, 1), \ \mbox{and therefore} \ m'(dx,0) = [1-p(x)]m(dx,1) + m(dx, 0), 
\eea
Roughly speaking, $m' \preceq m$ if the distribution $m'$ can be obtained by immediately stopping particles from the distribution $m$. In our context, $m_{t^-} = \dbP_{(X_t,I_{t-})}$ and $m_t = \dbP_{(X_t, I_t)}$, with $\dbP \in \cP(t,m)$, so that $m_t \preceq m_{t^-}$ with conditional transition probability $p(x) = p_t(x) := \dbP(I_t = 1 \mid X_t = x, I_{t-} = 1)$. When $p=\1_A$, with $A \in \cB(\dbR^d)$, we denote $m^A$ the probability measure defined by \eqref{order}.
We recall from \cite{TTZ2} that the set $\{m' : m' \preceq m\}$ is compact. Moreover, we have the following simple lemma whose proof is postponed to Appendix.

 \begin{lem}\label{pure-mixed}
 Let $m \in \cP_2(\bS)$ such that $x\in\dbR^d \to m\big((-\infty, x_1]\times\cds\times (-\infty, x_d], ~1\big)$ is continuous. Then $\{m^A : A \in \cB(\dbR^d)\}$ is $\cW_2$-dense in $\{ m' : m' \preceq m\}$.
 \end{lem}

In \cite{TTZ, TTZ2}, we see that $V$ is associated with the obstacle problem on Wasserstein space:
\bea\label{obstacle}
\min\Big\{\min_{m' \in C_U(t,m)} \big[-(\dbL U + F)(t,m')\big],~ (\dbD_I U)_* (t,m) \Big\} = 0, ~ (t,m) \in\mathbf{Q}_0,
\eea
with boundary condition $U|_{t=T} = g$. The function $(\dbD_IU)_*$ is the LSC envelope of 
 $$
 \dbD_IU:(t,m) \longmapsto \inf_{x\in\mbox{\tiny Supp}(m(.,1))}D_I U(t,m,x),
 $$
which is upper semicontinuous, but may not be continuous in general,  and the set 
 \beaa
 C_U(t,m) 
 &:=& 
 \big\{m' \preceq m : U(t,m') \ge U(t,m)\big\},
 ~~(t,m) \in \mathbf{Q}_0,
 \eeaa
indicates the set of positions at $t$ which improve $U$ by stopping the corresponding particles.

\subsection{Viscosity solutions}

For $\d > 0$ and $(t,m) \in \mathbf{Q}_0$, we introduce the neighborhood
$$ \cN_\d (t,m) := \big\{(s, \tilde m) :  s \in [t, t+\d], \dbP \in \cP(t,m), \ \tilde m \in \{ \dbP_{Y_{s-}}, \dbP_{Y_s} \} \big\}. $$
Note that $\cN_\d(t,m)$ is compact, as the closure of a càdlàg $\cP_2(\bS)$-valued graph.
For a  locally bounded function $U: \overline{\mathbf{Q}}_0\longrightarrow\dbR$, we introduce its LSC and USC envelopes relatively to $\cP(t,m)$, $U_*$ and $U^*$ respectively:
\beaa
U_*(t,m) := \underset{(s, \tilde m) \to (t,m)}{\lim \inf} U(s,\tilde m), \q U^*(t, m) := \underset{(s,\tilde m) \to (t,m)}{\lim \sup} U(s,\tilde m), 
&\mbox{for all}&
(t,m) \in \overline{\mathbf{Q}}_0, 
\eeaa
where the limits are taken on all sequences $\{t_n,m_n\}_{n \ge 1}$ converging to $(t,m)$, with $(t_n,m_n) \in \cN_{T-t}(t, m)$ 
for all $n$. Moreover, the smooth functions are in the following sense.
\begin{defn}\label{C12S}
Let $C^{1,2}_2(\overline\bQ_t)$ denote the set of continuous functions $U: \overline{\mathbf{Q}}_t \to \dbR$ such that 
\\
$~\hspace{5mm}\bullet$ $\pa_t U, \d_m U, \pa_x \d_m U_1, \pa_{xx}^2\d_m U_1$ exist, and are continuous in all variables, 
\\
$~\hspace{5mm}\bullet$ $\pa_{xx}^2\d_m U_1$ is bounded in $x$, locally uniformly in $(s, m)$,
\\
$~\hspace{5mm}\bullet$ $\pa_{\m \m}^2 U_{1,1}$ exists and $(s,m,x) \mapsto \pa_{\m \m}^2 U_{1,1} U(s,m,x,x)$  is continuous and has quadratic growth in $x$, locally uniformly in $(s,m)$. 
\end{defn}

We then introduce the sets of test functions
\bea
\overline{\cA}U(t,m) &:=& \Big\{\f \in C^{1,2}_2(\overline{\mathbf{Q}}_t):  (\f-U_*)(t,m) = \max_{\cN_\d(t, m)}(\f-U_*) \ \mbox{for some $\d > 0$}\Big\} \nonumber, \\ 
\underline{\cA}U(t,m) &:=& \Big\{\f \in C^{1,2}_2(\overline{\mathbf{Q}}_t): (\f-U^*)(t,m) = \min_{\cN_\d(t, m)}(\f-U^*) \ \mbox{for some $\d > 0$}\Big\}. \nonumber
\eea

\begin{defn}
\label{defn-viscosity}
Let $U : \ol\bQ_0 \rightarrow \dbR$ be locally bounded.

\no{\rm (i)} $U$ is a viscosity supersolution of \reff{obstacle} if, for any $(t, m)\in\mathbf{Q}_0$ and $m' \preceq m$,
\bea
\label{super}
U_*(t,m) \ge U_*(t,m') ~
\mbox{and}~
-(\dbL \f + F) (t,m) \ge 0,
~\mbox{for all $\f \in \overline{\cA}U(t,m)$.} 
 \eea
 
 \no{\rm (ii)} $U$ is a viscosity subsolution of \reff{obstacle} if,  for any  $(t, m)\in\mathbf{Q}_0$ s.t.  $C_{U^*}(t, m)  = \{m\}$,
\bea
\label{sub}
\min\{ -(\dbL \f + F) , (\dbD_I \f)_*\}(t,m) \le 0,\q \mbox{for all $\f \in \underline{\cA}U(t,m).$}
 \eea

\no{\rm (iii)} $U$ is a viscosity solution of \eqref{obstacle} if it is a viscosity supersolution and subsolution.

\end{defn}

\begin{rem} 
\label{rem-C12}
{\rm
In the present paper, the set $C^{1,2}_2(\ol\bQ_t)$ has a different definition from the one we used in our previous work \cite{TTZ2}, as we additionally require the existence of the second order derivative $\pa_{\m \m}^2 U_{1,1}$. However, we can see that the viscosity theory implied by requiring such regularity for the test functions is equivalent to the one developed in \cite{TTZ2}; in particular, the test functions involved in the proof of our comparison result \cite[Theorem 3.11]{TTZ2} are infinitely differentiable, so that the validity of this result does not depend on the definition of $C^{1,2}_2(\ol\bQ_t)$ we set.  \qed
}
\end{rem}

Recall the assumptions for our comparison result \cite[Theorem 3.11]{TTZ2}. Together with \cite[Theorem 3.9]{TTZ2} and Remark \ref{rem-C12}, this shows that $V$ is the unique viscosity solution of \eqref{obstacle}.
\begin{assum}\label{assum-comparison}
{\rm (i)} $b$ and $\si$ can be extended to $\cP_1(\bS)$ under $\cW_1$ continuously, and are uniformly Lipschitz continuous in $(x, m)$ under $\cW_1$;

\no{\rm (ii)} $f$ and $g$ extend continuously to $\cP_1(\bS)$ and have linear growth in $(x, m)$, under $\cW_1$;

\no{\rm (iii)}  $f$ is uniformly continuous in $(t,x,m)$, under $\cW_2$ for $m$; and all first and second order derivatives of $\si(t, x, m)$ with respect to $x$ and those of $\d_m \si(t, x, m, \tilde x, \tilde i)$ with respect to $(x, \tilde x)$ are bounded, continuous in $t$, and Lipschitz continuous in $(x, \tilde x, m)$, under $\cW_1$ for $m$. 
\end{assum}

\section{The $N$-player optimal stopping problem}
\label{sect-N}
In this section we introduce a general $N$-player optimal stopping problem, following Kobylanski, Quenez \& Rouy-Mironescu \cite{KQR}. In the symmetric case, this reduces to the $N$-player problem corresponding to the mean field problem \reff{weakoptstop}, as we will see in Subsection \ref{sect:sym}.

\subsection{The general $N$-player optimal stopping problem}
\label{sect-general}
Let $(\O^0, \cF_T^0, \dbF^0, \dbP^0)$ be a filtered probability space, equipped with a sequence of $d$-dimensional Brownian motions $\{W^k\}_{k \ge 1}$. For each $N\ge 1$, let $\dbF^N$ denote the canonical filtration of $\{W^k\}_{1 \le k \le N}$,  $\cT_{t,T}^N$ the set of  $[t,T]-$valued $\dbF^N-$stopping times for all $t \in [0,T]$,  $\mathbfcal{T}_{t,T}^N := (\cT_{t,T}^N)^N$, $[N] := \{1, \dots, N\}$. We emphasize that, unlike the weak formulation in the previous section, here we use strong formulation. In particular, the stopping times in $\cT_{t,T}^N$ are adapted to the Brownian motion filtration, in the spirit of open loop controls. Moreover, in this section $N$ is fixed, but in the next section we will send $N$ to $\infty$.

Let $(b_k,\sigma_k, f_k):[0,T]\times \bS^N  \longrightarrow\dbR^d\times\cS_d^+(\dbR)\times \dbR$, $k \in [N]$,  and $g_N: \dbR^{d \times N}\to \dbR$ satisfy:
\begin{assum}
\label{assum-N}
\no{\rm (i)} The functions $b_k, \si_k$, $k\in [N]$,  are continuous in $t$ and uniformly Lipschitz continuous in $\bx$ (with the Lipschitz constant possibly depending on $N$);

\no{\rm (ii)} The functions $f_k$, $k\in [N]$,  and  $g_N$ are continuous in $t$ and uniformly continuous in $\bx$ (with the uniform continuity possibly depending on $N$). 
\end{assum}

 Fix $t\in [0, T)$ and $\by = (\bx, \bi) \in \bS^N$.  For  $\btau := (\t_1, \dots, \t_N) \in \mathbfcal{T}_{t,T}^N$,  let $\bX^\btau = \bX^{N, \btau}  = (X^{\btau, 1}, \dots, X^{\btau, N})$ and $\bI^\btau = \bI^{N, \btau} = (I^{\btau, 1}, \dots, I^{\btau, N})$ solve the following system of SDEs: 
\bea\label{Ndynamics}
X_s^{\btau, k} = x_k + \int_t^s b_k(r, \bY_r^{\btau}) I_r^{\btau, k} dr + \int_t^s \si_k(r, \bY_r^{\btau}) I_r^{\btau, k} dW_r^k, ~ I_s^{\btau, k} = i_k \1_{\t_k > s},  \q \dbP^0\mbox{-a.s.} 
\eea
for $k\in [N]$, where $\bY^\btau = \bY^{N, \btau} := (\bX^\btau, \bI^\btau)$. Here and throughout the paper, the dependence of $\bY$ on $N$ is omitted for notational simplicity. Moreover, in light of \reff{barP}, let $\bar \bY=(\bar\bX, \bar \bI)$ denote the unstopped version, namely $\tau_k = T$, or equivalently $\bar I^k_s = i_k$, for all $k\in [N]$.

We now consider the problem
\bea\label{Nparticles}
\left.\ba{c}
\dis v^N(t, \by) :=  \sup_{\btau \in \mathbfcal{T}_{t,T}^N}  J_N(t,\by, \btau),\q (t, \by)\in \ol\L^N_0,\\
\dis \mbox{where}\q J_N(t,\by, \btau):= 
  \dbE^{\dbP^0}_{t, \by}\Big[ \frac 1 N \sum_{k=1}^N \int_t^T f_k(r, \bY_r^\btau) I_r^{\btau, k}dr + g_N\big(\bY_T^\btau \big) \Big].
\ea\right.
\eea
Here the subscript $(t, \by)$ in $ \dbE^{\dbP^0}_{t,\by}[\cdot] :=  \dbE^{\dbP^0}_{t,\by}[\cdot | \bY_{t-}^\btau = \by]$ indicates the initial condition. 

\subsection{Viscosity solutions with tangency in mean}
\label{sec:viscosityN}
The function $v^N$ is associated with the following partial differential equation on $\ol\L_0^N$:
\bea\label{cascade}
\min\{-\pa_t u(t,\by) - \bi \cdot (\cL u + \bff)(t,\by), u(t,\by) - \max_{\bi' < \bi} u(t,\bx,\bi') \} = 0,  
\eea
with \by = (\bx, \bi) and boundary condition $u |_{t=T} = u |_{\bi=0} = g_N$, where  $\bi' < \bi$ means $i_k' \le i_k$, $k \in [N]$, and $\bi \neq \bi'$, $\bff := (f_1,\cds, f_N)$, and
\bea\label{Ndiffop}
\cL \f := \Big(b_k \cd \pa_{x_k}\f + \frac 1 2 \si_k^2: \pa_{x_k x_k}^2\f \Big)_{k \in [N]} \q \mbox{for all $\f \in C^2(\dbR^{d \times N})$}.
\eea
We call \eqref{cascade} \textit{cascade obstacle equation}. We consider the notion of viscosity solutions introduced by Ekren, Keller, Touzi \& Zhang \cite{EKTZ} where, for any function $u: \ol \L_0^N \to \dbR$, the corresponding set of test functions at some point $(t,\by)\in \ol \L_0^N$ is defined by a tangency mean rather than the standard pointwise tangency of Crandall \& Lions:
{\small
\bea\label{Ntest}
\left.\ba{c}
\dis \overline{\cA}u(t,\by):=  \Big\{\phi \in C^{1,2}(\ol\L_t^N):  \exists \ch\in \cT_{t,T}^+, \ (\phi-u)(t,\by) =  \max_{\th \in \cT_{t,T}^N} \dbE^{\dbP^0}_{t,\by}\Big[ (\phi - u)(\th\wedge \ch, \bar \bY_{\th\wedge \ch}) \Big]\Big\}, \\ 
\dis \underline{\cA}u(t,\by) := \Big\{\phi \in C^{1,2}(\ol\L_t^N):  \exists \ch \in \cT_{t,T}^+, \ (\phi-u)(t,\by) =  \min_{\th \in \cT^N_{t,T}} \dbE^{\dbP^0}_{t, \by}\Big[ (\phi - u)(\th\wedge\ch, \bar \bY_{\th\wedge\ch}) \Big]\Big\},
\ea\right.
\eea}
with $\cT_{t,T}^+ := \{ \ch \in \cT^N_{t,T} : \ch > t \}$, and $ C^{1,2}(\ol\L_t^N)$ denoting the set of continuous test functions $\phi$ such that: for any $\bi =(i_1, \cds, i_N)$ and for all $k\in [N]$ with $i_k =1$, the derivatives $\pa_t \phi(s, \bx, \bi)$, $\pa_{x_k} \phi(s, \bx, \bi), \pa_{x_k, x_k} \phi(s, \bx, \bi)$ exist and are continuous for $(s, \bx) \in [t, T]\times \dbR^{d\times N}$.

\begin{defn}\label{defn:viscosityN}
Let $u : \ol\L_0^N \longrightarrow \dbR$. \\
{\rm(i)} $u$ is a viscosity supersolution of \eqref{cascade} if, for all $(t,\by) \in \ol\L_0^N$,
$$ \min\{-\pa_t \phi(t,\by) - (\cL \phi + \bff)(t,\by) \cdot \bi, u(t,\by) - \max_{\bi' < \bi} u(t,\bx,\bi') \} \ge 0 \q \mbox{for all $\phi \in \ol \cA u(t,\by)$.} $$
{\rm(ii)} $u$ is a viscosity subsolution of \eqref{cascade} if, for all $(t,\by) \in \L_0^N$,
$$ \min\{-\pa_t \phi(t,\by) - (\cL \phi + \bff)(t,\by) \cdot \bi, u(t,\by) - \max_{\bi' < \bi} u(t,\bx,\bi') \} \le 0 \q \mbox{for all $\phi \in \ul \cA u(t,\by)$.} $$
{\rm (iii)} $u$ is a viscosity solution of \eqref{cascade} if it is a viscosity supersolution and subsolution.
\end{defn}
We then have the following result, whose proof is postponed  to \S \ref{sect-Ndimpb}.
\begin{thm}
\label{thm:Nexistence}
Under Assumption \ref{assum-N}, the function $v^N$ in \reff{Nparticles} is the unique continuous viscosity solution of \eqref{cascade} with boundary conditions $v^N|_{t=T} = v^N|_{\bi = 0} = g$.
\end{thm}

\begin{rem}
\label{rem-viscosityN}
{\rm We note that a viscosity solution in our sense is always a viscosity solution in Crandall-Lions sense, as the latter involves a smaller class of test functions. In particular, the comparison principle and the uniqueness of viscosity solution in the latter sense implies the same results in our sense. Our definition turns out to be more convenient for our proof of the convergence results in the next section. 
\qed}
\end{rem}

\subsection{The symmetric case}
\label{sect:sym}
Recall the functions $b, \si, f, g$ in \S \ref{sect-MF} and the notation $m^N$ in \reff{QN}. In this subsection we consider the symmetric case:
\bea
\label{symmetric}
\f_k(t, \by) = \f(t, x_k, m^N(\by)),~\mbox{for}~ \f = b, \si, f,\q\mbox{and}\q g_N(\by) = g(m^N(\by)).
\eea
We first have the following simple lemma, whose  proof is postponed to \S \ref{sect-Ndimpb}.

\begin{lem}\label{lem:unif-estimate}
{\rm (i)} Under Assumption \ref{assum-bsig} \rm{(i)} and \reff{symmetric}, there is a constant $C \ge 0$, independent of $N$, such that, for $(t,\by)\in \L^N_0$, $\btau \in\mathbfcal{T}_{t,T}^N$, $k\in [N]$, and $\th_1, \th_2\in \cT^N_{[t, T]}$ with $\th_1 \le \th_2$,
\bea
\label{unif-estimate11}
&\dis \dbE\Big[ \sup_{t\le s\le T}\lvert X_s^{\btau, k} \rvert^2 \Big] \le C(1+|x_k|^2 + \|\bx\|_2^2);\ms\\ 
\label{unif-estimate12}
&\dis \dbE\Big[ \sup_{\th_1 \le s\le \th_2}\lvert X_s^{\btau, k} - X_{\th_1}^{\btau, k}  \rvert^2 \Big] \le C(1+|x_k|^2 + \|\bx\|_2^2) \|\th_2-\th_1\|_\infty.
 \eea
 
\no{\rm (ii)} If we assume further that  Assumption \ref{assum-comparison} (i) holds, then
 \bea
\label{unif-estimate21}
&\dis \dbE\Big[ \sup_{t\le s\le T}\lvert X_s^{\btau, k} \rvert \Big] \le C(1+ |x_k|+ \|\bx\|_1);\ms\\ 
\label{unif-estimate22}
&\dis \dbE\Big[ \sup_{\th_1 \le s\le \th_2}\lvert X_s^{\btau, k} - X_{\th_1}^{\btau, k}  \rvert \Big] \le C(1+|x_k|+\|\bx\|_1) \sqrt{\|\th_2-\th_1\|_\infty}.
 \eea
\end{lem}

Next, under \reff{symmetric} one can easily see that $v^N$ is also symmetric in $\by$, that is, there exists a function $V^N: \ol \bQ^N_0\to \dbR$ such that
\bea
\label{VNsymmetric}
v^N(t, \by) = V^N(t,  m^N(\by)).
\eea

\begin{defn}
\label{defn-viscositysymmetric}
Under \reff{symmetric}, we call a function $U^N: \ol \bQ^N_0\to \dbR$ a symmetric viscosity solution (resp. subsolution, supersolution) of \reff{cascade} if $u^N(t, \by):= U^N(t, m^N(\by))$ is a viscosity solution (resp. subsolution, supersolution) of \reff{cascade} in the sense of Definition \ref{defn:viscosityN}.
\end{defn}
 The following result establishes stronger regularity than Theorem \ref{thm:Nexistence} and its proof is also postponed to  \S \ref{sect-Ndimpb}.

\begin{thm}
\label{thm:Nexistence-symmetric}
 Let Assumptions \ref{assum-bsig}, \ref{assum-comparison} {\rm(i)-(ii)} and \reff{symmetric} hold. 
 
\no{\rm (i)} The function $v^N$ is the unique continuous symmetric viscosity solution of \eqref{cascade}. 

\no{\rm (ii)} There exist a constant $C$,  independent of $N$,  such that 
 \bea
 \label{unif-estimate31}
  |V^N(t, m_N)| \le C\big(1+ \|m_N\|_1\big).
 \eea
 
 \no{\rm (iii)} For any $R>0$, there exist a modulus of continuity function $\rho_R$, which may depend on $R$ but not on $N$,  such that, for any $(t, m_N), (\tilde t, \tilde m_N)\in  \ol \bQ^N_0$ with $\|m_N\|_2, \|\tilde m_N\|_2 \le R$,
 \bea
\label{unif-estimate32}
 |V^N(t, m_N)- V^N(\tilde t, \tilde m_N)| \le \rho_R\big(\cW_2(m_N, \tilde m_N) +|\tilde t-t|\big).
 \eea
\end{thm}

Moreover, similar to Carmona \& Delarue \cite[Vol. 1, Propositions 5.35 \& 5.91]{CarDel}  we have the following result, and we shall provide a proof in Appendix.
\begin{lem}
\label{lem-derivative}
Assume $\f\in C^{1,2}_2(\ol \bQ_0)$ and $\phi(t, \by) := \f\big(t, m^N(\by)\big)$ for all $(t, \by) = (t, \bx, \bi) \in [0,T] \times \bS^N$. Then, for all $(s, \by) \in \ol\L^N_0$ and $k \in [N]$ such that $i_k = 1$,
\bea\label{connection-derivatives}
\left.\ba{c}
\dis \pa_{x_k}\phi(s, \by) = \frac 1 N \pa_x \d_m \f_1(s, m^N(\by), x_k), \\ 
\dis \pa_{x_k x_k}^2 \phi(s, \by)  = \frac 1 N \pa_{xx}^2 \d_m \f_1(s, m^N(\by), x_k) + \frac{1}{N^2} \pa_{\m \m}^2 \f_{1,1}(s, m^N(\by), x_k, x_k). \ea\right.
\eea
\end{lem}

\section{Main results}
\label{sect-main}
Our first main result states the convergence of the value function of the multiple optimal stopping problem to the corresponding mean field problem. 

 \begin{thm}\label{thm-vfconv}
Let Assumptions \ref{assum-bsig} and \ref{assum-comparison} hold, and $V^N$ and $V$ be respectively the value functions of the multiple and mean field optimal stopping problems \eqref{Nparticles}-\reff{symmetric}-\reff{VNsymmetric} and \eqref{weakoptstop}. Then $V^N$ converges to $V$, i.e.
\bea
\label{vfconv}
\underset{\tiny \begin{array}{c} N \to \infty, s \to t \\ m_N\in \cP^N(\bS) \ \overset{\cW_2}{\longrightarrow} m \end{array}}{\lim} V^N\big(s, m_N\big) = V(t,m) ,\q \mbox{for all $(t,m) \in \ol \bQ_0$.}
\eea
\end{thm}

The above theorem relies heavily on the following result. For a sequence of functions $U^N:  \ol \bQ^N_0\to \dbR$ which are locally bounded functions, uniformly in $N$, we denote 
\bea
\label{wbar}
 \ul U(t,m) :=  \!\!\!\!\!\underset{\tiny \begin{array}{c} N \to \infty, s \to t \\ m_N\in \cP^N(\bS) \overset{\cW_2}{\longrightarrow} m \end{array}}{\lim \inf}   \!\!\!\!\! U^N(s, m_N), \q \mbox{and} \q \ol U(t,m) := \!\!\!\!\! \underset{\tiny \begin{array}{c} N \to \infty, s \to t \\ m_N\in \cP^N(\bS) \overset{\cW_2}{\longrightarrow} m\end{array}}{\lim \sup} \!\!\!\!\!   U^N(s, m_N),
 \eea
for all $(t,m) \in \ol\bQ_0$.  We then have

\begin{thm}\label{thm-visc-cv}
  Let Assumption \ref{assum-bsig} and \eqref{symmetric} hold and $U^N: \ol \bQ^N_0\to \dbR$  be a sequence of continuous and locally bounded functions, uniformly in $N$. 

\no{\rm (i)} If each $U^N$  is a viscosity subsolution of \eqref{cascade}, then $\ol U$ is a viscosity subsolution of $\eqref{obstacle}$. 

\no{\rm (ii)} If each $U^N$  is a viscosity supersolution of \eqref{cascade}, and $\ul U$ is continuous, then $\ul U$ is a viscosity supersolution of $\eqref{obstacle}$. 
\end{thm}
The above two theorems are proved in \S \ref{sect-convergence}. We now turn to the convergence of the optimal strategies. 
Our second main result is as follows.

\begin{thm}\label{thm:propagation-optimal}
Let Assumptions \ref{assum-bsig}  and \ref{assum-comparison}  hold.  Fix $(t, m)\in \bQ_0$ and  $\by^N\in \bS^N$, $N \ge 1$, such that $m^N(\by^N) \overset{\cW_2}{\longrightarrow} m$. Recall \reff{Ndynamics} and assume  \reff{symmetric} holds.

\no {\rm (i)} Denote $m^N(\bY):= {1\over N}\sum_{k=1}^N \d_{Y^k}\in \cP_2(\O)$. Assume $\hat \btau^N \in \mathbfcal{T}_{t,T}^N$ is $\e_N$-optimal for \eqref{Nparticles} for all $N \ge 1$, with $\e_N \to 0$. Then, $\big\{\dbP^0 \circ \big(m^N(\bY^{\hat \btau^N})\big)^{-1} \big\}_{N \ge 1}$ is tight, and all its accumulation points are supported on the set of the optimal controls for \eqref{weakoptstop}. \\ 
{\rm (ii)} Let $\dbP^* \in \cP(t,m)$ be optimal for \eqref{weakoptstop}.  Then there exist $\e_N \to 0$ and $\hat \btau^N \in \mathbfcal{T}_{t,T}^N$ which is $\e_N$-optimal for \eqref{Nparticles} such that 
$\dis m^N(\bY^{\hat \btau^N})  \overset{\cW_2}{\longrightarrow} \dbP^*$, $\dbP^0$-a.s.\\
{\rm (iii)}  Assume further that \eqref{weakoptstop} has a unique optimal control $\dbP^*$. Then for any sequence of $\e_N$-optimal stopping strategies $\{\hat \btau^N\}_{N \ge 1}$  for \eqref{Nparticles},  $\dis m^N(\bY^{\hat \btau^N})  \overset{\cW_2}{\longrightarrow} \dbP^*$, $\dbP^0$-a.s.
 \end{thm}
Since $N$ varies here, we use the superscript $^N$ in $\hat \btau^N$ to indicate its dependence on $N$.

%

Theorem \ref{thm:propagation-optimal}  relies heavily on the following result concerning the propagation of chaos for arbitrary stopping strategies. Both results are proved in \S  \ref{sec-proof-propagation}.
 
\begin{thm}\label{thm:approx}
Let Assumption \ref{assum-bsig}  hold, and fix $(t, m)\in \bQ_0$ and  $\by^N\in \bS^N$, $N \ge 1$, such that $m^N(\by^N) \overset{\cW_2}{\longrightarrow} m$. Recall \reff{Ndynamics} and assume  \reff{symmetric} holds.

\no{\rm (i)}  For any $\btau^N \in \mathbfcal{T}_{t, T}^N$, the sequence $\big\{\dbP^0 \circ \big(m^N(\bY^{\btau^N})\big)^{-1} \big\}_{N \ge 1}$ is tight, and all its accumulation points  are supported on $\cP(t,m)$. \\
{\rm (ii)} Assume further that $b, \si$ are uniformly Lipschitz continuous in $m$ under $\cW_1$. Then, for any $\dbP \in \cP(t,m)$, there exists $\btau^N \in \mathbfcal{T}_{t,T}^N$, $N \ge 1$, such that $m^N(\bY^{\btau^N})  \overset{\cW_2}{\longrightarrow} \dbP$, $\dbP^0$-a.s.
\end{thm}

\section{Convergence of the value function}
\label{sect-convergence}
In this section we prove Theorems \ref{thm-vfconv} and \ref{thm-visc-cv}.

\ms
\no {\bf Proof of Theorem \ref{thm-vfconv}} (given Theorem \ref{thm-visc-cv}). First, for  any $(t, m)\in \ol \bQ_0$ and $(t^N, m_N(\by^N)) \to (t, m)$ as $N\to \infty$, the sequence $\{V^N(t^N, m_N(\by^N))\}_{N\ge 1}$ is bounded by \reff{unif-estimate31}. Then the corresponding functions $\ol V, \ul V$ defined as in \reff{wbar} are finite, and are continuous by Theorem \ref{thm:Nexistence-symmetric} (iii). In particular, $\ol V(T,\cdot) = \ul V(T,\cdot) = g$. Applying Theorem \ref{thm-visc-cv} we see that $\ol V$ and $\ul V$ are a viscosity subsolution and supersolution of  \eqref{obstacle}. Then it follows from  the comparison principle \cite[Theorem 3.11]{TTZ2}  that $\ul V \ge \ol V$. As the converse inequality holds by definition, this shows that $\ul V = \ol V$, and by uniqueness, they are equal to $V$ and thus \reff{vfconv} holds. \qed

\ms
\no{\bf Proof of Theorem \ref{thm-visc-cv}.} Without loss of generality, we prove the theorem only at $t=0$. 

(i) 
Fix $m \in \cP_2(\bS)$ and $\f \in \ul \cA \ol U(0,m)$ with corresponding $\d_0 \in (0,T)$. We assume w.l.o.g. (see \cite[Remark 3.2]{TTZ2}) that 
\bea\label{strict-min}
\mbox{$(0,m)$ is a strict minimum of $(\f - \ol U)$ on $\cN_{\d_0}(0,m)$, with $C_{\ol U}(0,m) = \{m\}$. }
\eea
Let $(t^N, \by^N)$ be such that 
\bea
\label{mNconv}
(t^N, m^N(\by^N)) \longrightarrow (0,m),\q\mbox{and}\q U^N(t^N, m^N(\by^N)) \longrightarrow \ol U(0,m),\q\mbox{as}\q N \to \infty.
\eea
 We also introduce  the functions $\phi^N(s,\by) := \f(s, m^N(\by))$ for all $N\ge 1$ and $(s,\by) \in \ol \L^N_{t^N}$. 

\ms
\no{\bf Step 1:} We prove that, without loss of generality, $\{\by^N = (\bx^N, \bi^N)\}_{N \ge 1}$ may be taken s.t.
\bea\label{strict-positive}
U^N(t^N, \by^N) - \max_{\bi' < \bi^N} U^N(t^N, \bx^N, \bi') > 0 \q \mbox{for all $N \ge 1$.}
\eea
Indeed, for each $N\ge 1$, if $\bi^N = 0$, then the maximum is equal to $-\infty$ and thus \eqref{strict-positive} is trivially satisfied. Now assume $\bi^N \neq 0$. Then there exists $\bi^{N,*} \le \bi^N$ attaining the maximum $ \max_{\bi' \le \bi^N} U^N(t^N, \bx^N, \bi')$. Note that $\bi^{N,*}$ can be chosen to be minimal in $\mathrm{arg max}_{\bi' \in \bi^N} U^N(t^N, \bx^N, \bi')$ for the partial order on $\{0,1\}^N$. This implies that $\by^{N,*} := (\bx^N, \bi^{N,*})$ satisfies \reff{strict-positive}. 

It remains to show that $\by^{N,*}$ also satisfies \reff{mNconv}. First, as $m^N(\by^N)$ converges, the first marginal of $m^N(\by^{N,*})$ converges. Next, since its second marginal is a measure on $\{0,1\}$, the sequence $\{m^N(\by^{N,*})\}_{N \ge 1}$ is tight, and we may thus extract a subsequence (still denoted the same) converging to some $m^* \preceq m$. Then, noting that $U^N(t^N, \by^N)  \le U^N(t^N, \by^{N,*})$, by  taking the $\underset{N \to \infty}{\lim \inf}$ we obtain
\beaa
\ol U(0,m) &=& \lim_{N \to \infty} U^N(t^N, m^N(\by^N)) \le \underset{N \to \infty}{\lim \inf} \ U^N(t^N, m(\by^{N,*})) \\
&\le& \underset{N \to \infty}{\lim \sup} \ U^N(t^N, m(\by^{N,*})) \le \ol U(0,m^*), 
\eeaa
Then $m^* \in C_{\ol U}(0,m)$, and therefore $m^* = m$ by \eqref{strict-min}, so that equality holds everywhere and and $ U^N(t^N, m^N(\by^{N,*})) \underset{N \to \infty}{\longrightarrow} \ol U(0,m)$. Consequently, $\{\by^{N,*}\}_{N \ge 1}$ satisfies \reff{mNconv}, and we may replace $\{\by^N\}_{N \ge 1}$ with $\{\by^{N,*}\}_{N \ge 1}$, which satisfies in addition \eqref{strict-positive}.

\ms
\no{\bf Step 2:} Fix $\d \in (0, {\d_0\over 3})$. By standard optimal stopping theory, there exists $\th_\d^N \in \cT_{t^N,T}$ s.t.
$$ 
\dbE_{t^N, \by^N}^{\dbP^0}\Big[(\phi^N - U^N)(\th^N_\d \wedge \ch_\d^N, \bar \bY^N_{\th_\d^N \wedge \ch_\d^N}) \Big] = \inf_{\th \in \cT_{t^N, T}} \dbE_{t^N, \by^N}^{\dbP^0}\Big[ (\phi^N - U^N)(\th \wedge \ch_\d^N, \bar \bY^N_{\th \wedge \ch_\d^N}) \Big], 
$$
where $\bar \bY^N := (\bar \bX^N, \bar \bi^N)$ is the unstopped version of \eqref{Ndynamics} starting from $(t^N, \by^N)$, and
\begin{align*}
 \ch_\d^N := \inf \Big\{ s \ge t^N : \cW_2\big(m^N(\bar \bY_s), m^N(\by^N) \big) = 2\d  \Big\} \wedge (t^N + 2\d).
\end{align*}
Since $\{U^N\}_{N \ge 1}$ is locally bounded, we may assume w.l.o.g. that $N$ is large enough and $\d$ is small enough so that $\{(\phi^N - U^N)(\th_\d^N \wedge \ch_\d^N, \bar \bY^N_{\th_\d^N \wedge \ch_\d^N}) \}_{N \ge 1}$ is uniformly bounded.

\ms
\no{\bf Step 3:} We next show that, for any $\d$, 
\bea\label{theta}
\underset{N \to \infty}{\lim \sup} \ \dbP^0\Big(\th_\d^{N} < \ch_\d^N \Big) > 0.
\eea
Indeed, assume to the contrary that $\underset{N \to \infty}{\lim \sup} \ \dbP^0\Big(\th_\d^{N} < \ch_\d^N \Big)=0$. Then
\bea\label{estsup}
(\f - \ol U)(0,m) &=& \underset{N \to \infty}{\lim} (\f - U^N)(t^N, m^N(\by^N)) \nonumber \\
 &\ge& \underset{N \to \infty}{\lim \sup} \ \dbE_{t^N, \by^N}^{\dbP^0}\Big[(\f - U^N)(\th_\d^N \wedge \ch_\d^N, m^N(\bar \bY^N_{ \th_\d^N \wedge \ch_\d^N}))\Big] \nonumber  \\
&\ge& \underset{N \to \infty}{\lim \inf} \ \dbE_{t^N, \by^N}^{\dbP^0}\Big[(\f - U^N)(\th_\d^N \wedge \ch_\d^N, m^N(\bar \bY^N_{\th_\d^N \wedge \ch_\d^N}))\Big] \nonumber  \\
&=&\underset{N \to \infty}{\lim \inf} \ \dbE_{t^N, \by^N}^{\dbP^0}\Big[ \Big\{ (\f - U^N)(\ch_\d^N, m^N(\bar \bY^N_{\ch_\d^N}))(1-\1_{\{\th_\d^N < \ch_\d^N\}})  \\
&&\hspace{2cm}+ (\f - U^N)( \th_\d^N , m^N(\bar \bY^N_{\th_\d^N }))\1_{\{\th_\d^N < \ch_\d^N\}} \Big\} \Big] \nonumber \\
&=& \underset{N \to \infty}{\lim \inf} \  \dbE_{t^N, \by^N}^{\dbP^0}\Big[(\f - {\rm U^N})(\ch_\d^N,  m^N(\bar \bY^N_{\ch_\d^N}) )\Big], \nonumber
\eea
where we use the uniform boundedness of $\{(\phi^N - U^N)(\th \wedge \ch_\d^N, \bar \bY^N_{\th \wedge \ch_\d^N}) \}_{N \ge 1}$.

Recall that $ m^N(\bar \bY^N) \in \cP_2(\O)$, $\dbP^0$-a.s. 
Since $(t^N, m^N(\by^N)) \underset{N \to \infty}{\longrightarrow} (0,m)$ and $\bar \bY^N$ is an unstopped diffusion whose coefficients satisfy the usual Lipschitz conditions, by classical propagation of chaos (see Oelschläger \cite{Oelschlager}) and compactness of $[0,T]$, again after possibly passing to a subsequence, we see that $\n^N := \dbP^0 \circ \big( \ch_\d^N, m^N(\bar \bY^N) \big)^{-1}$ converges weakly to some element of $\n \in \cP_2( [0,T] \times  \cP_2(\O))$. Let $(\t, \l)$ denote the canonical mapping on $[0,T] \times \cP_2(\O)$. Recall that $Y$ is the canonical mapping on $\O$, and since $\bar \bY^N$ is unstopped, then $Y$ is unstopped, $m^N(\bar \bY^N)(\o^0)$-a.s., for $\dbP^0$-a.e. $\o^0\in \O^0$. This implies that  $Y$ is unstopped, $\l(\ol \o)$-a.s., for $\n$-a.e. $\ol \o \in  \cP_2( [0,T] \times  \cP_2(\O))$, and thus the mapping $(\t, \l) \mapsto \l_{Y_\t}$ (the $\l$-distribution of the random variable $Y_\t$)  is continuous, $\nu$-a.s. Then, by the definition of $\overline{U}$, together with Fatou's lemma, we deduce from \reff{estsup} that
\bea
\label{estsup2}
(\f - \ol U)(0,m)\ge \underset{N \to \infty}{\lim \inf} \ \dbE^{\n^N}\Big[ (\f - U^N)(\t, \l_{Y_\t})\Big] \ge\dbE^{\n}\Big[ (\f - \ol U)(\t, \l_{Y_\t})\Big]. 
\eea
Thus  there exists  $\bar \o \in [0,T] \times \cP_2(\O)$ such that  $Y$ is  unstopped under $\l(\bar\o)$, and
\beaa
(\f - \ol U)(0,m)
&\ge& 
(\f - \ol U)(\t(\bar \o), [\l(\bar\o)]_{Y_{\t(\bar \o)}}).
\eeaa
Moreover, we may choose $\bar\o$ so that
\beaa
\cW_2\big([\l(\bar\o)]_{Y_{\t(\bar \o)}}, m \big) \vee \t(\bar\o) \le 3\d <\d_0,\q \cW_2\big([\l(\bar\o)]_{Y_{\t(\bar \o)}}, m \big) + \t(\bar\o) \ge \d,
\eeaa
as we have  by the definition of $\ch_\d^N$, for all $N\ge 1$ that
\beaa
\cW_2\big(\l_{Y_\t}, m^N(\by^N) \big)\vee (\t-t^N) \le 2\d \le \cW_2\big(\l_{Y_\t}, m^N(\by^N) \big)+ (\t-t^N), \q \mbox{$\n^N$-a.s.,}
\eeaa
and for sufficiently large $N$:
\bea\label{lastequation}
\cW_2\big(\l_{Y_\t}, m \big)\vee \t \le 3\d 
&\mbox{and}&
\d\le \cW_2\big(\l_{Y_\t}, m \big)+ \t, \q \n^N-\mbox{a.s.}
\eea
By the continuity of $\cW_2$ and $\bar\o\longmapsto\l_{Y_\t}(\bar\o)$, the last inequalities exhibit a fixed closed support for $\nu^N$, which is then inherited by the weak limit $\n$. Note further $\l_{Y_{t-}} = m$, we deduce that $(0, m) \neq (\t(\bar \o),[\l(\bar\o)]_{Y_{\t(\bar \o)}}) \in \cN_{\d_0}(0,m)$. This 
contradicts with \reff{strict-min} that $(0, m)$ is a strict minimum on $\cN_{\d_0}(0,m)$. Thus, \eqref{theta} holds true.

\ms
\no{\bf Step 4:} We prove that $\phi^N$ is a test function of $U^N$ at some point. Introduce
$$\cZ_s^N := \essinf_{\th \in \cT_{t^N, T}^N} \dbE^{\dbP^0}\big[ (\phi^N - U^N)(s, \bar \bY_s^{t^N, \by^N}) | \cF_s^N \big] \ \mbox{for all $s \in [t^N, T]$,}$$ 
i.e., the $\dbP^0$-lower Snell envelope of $s \mapsto (\phi^N - U^N)(s, \bar \bY_s^{t^N, \by^N})$. Note that $\cZ^N$ is a $\dbP^0$-submartingale. By \eqref{theta}, for all $\d > 0$ and $N \ge 1$, we may find $\o^{\d,N} \in \O^0$ such that 
$\th_\d^N(\o^{\d,N}) < \ch_\d^N(\o^{\d,N})$ and, denoting $(t_\d^N, \by_\d^N) := (\th^{\d,N}(\o^{\d,N}), \bar \bY_{ \th_\d^N(\o^{\d,N})}^{t^N,\by^N}(\o^{\d,N}))$,
\beaa
(\phi^N - U^N)(t_\d^N, \by_\d^N) &=& \cZ_{t_\d^N}^N(\o^{\d,N}) \le \dbE^{\dbP^0}\big[ \cZ_{\th  \wedge \ch_\d^N}^N \big] \le \dbE^{\dbP^0}\big[(\phi^N - U^N)(\th  \wedge \ch_\d^N, \bar \bY_{\th  \wedge \ch_\d^N}^{t^N, y^N}) \Big],
\eeaa
for all $\th \in \cT_{t_\d^N,T}$. Thus we have
\bea\label{sub-visc-tg} 
(\phi^N - U^N)(t_\d^N, \by_\d^N) = \min_{\tiny \begin{array}{c} \th \in \cT_{t_\d^N, T} \end{array}} \dbE_{t_\d^N, \by_\d^N}^{\dbP^0}\Big[(\phi^N - U^N)(\th \wedge \ch_\d^N, \bar \bY_{\th \wedge \ch_\d^N}) \Big]. 
\eea
Therefore, since $\phi^N \in C^{1,2}(\ol\L_{t^N}^N)$ by Lemma \ref{lem-derivative}, we have $\phi^N \in \ul \cA^N U^N(t_\d^N, \by_\d^N)$, and thus the subsolution property of $U^N$ provides, 
{\small
\bea\label{subcv}
\min\Big\{-\pa_t \phi^N(t_\d^N, \by_\d^N) - (\cL \phi^N + \bff^N)(t_\d^N, \by_\d^N) \cdot  \bi^N, U^N(t_\d^N, \by_\d^N) - \max_{\bi' < \bi^N} U^N(t_\d^N, \bx_\d^N, \bi')\Big\} \le 0.  
\eea
}
\no{\bf Step 5:} We finally derive the viscosity subsolution property of $\ol U$. First observe that, since $t_\d^N = \th_\d^N(\o^{\d,N}) < \ch_\d^N(\o^{\d,N})$, we have
$ (t_\d^N, \by_\d^N) \underset{\d \to 0}{\longrightarrow} (t^N,\by^N). $
Since $\phi^N$, $\cL \phi^N + \bff^N$ and $U^N$ are continuous, sending $\d$ to $0$ in \eqref{subcv} provides
$$ 
\min\Big\{\!\! -\pa_t \phi^N(t^N, \by^N) - (\cL \phi^N + \bff^N)(t^N, \by^N) \cdot  \bi^N\!\!, ~U^N(t^N, \by^N) - \max_{\bi' < \bi^N} U^N(t^N, \bx^N, \bi')\Big\} \le 0.  $$
By \reff{strict-positive}, this implies
\bea\label{subcv-phi}
-\pa_t \phi^N(t^N, \by^N) - (\cL \phi^N + \bff^N)(t^N, \by^N) \cdot  \bi^N \le 0, \q \mbox{for all $N \ge 1$.}
\eea
Thus, noting that $i^N_k$ takes only values $0$ and $1$, by \eqref{connection-derivatives} and \reff{symmetric} we have
\bea
\label{subcv-phi2} 
&&\dis -\pa_t \phi^N(t^N, \by^N) - (\cL \phi^N + \bff^N)(t^N, \by^N) \cdot  \bi^N   = -\dbL  \f(t^N, m^N(\by^N))\\
&&\dis - {1\over N} \sum_{k  = 1}^N f(t^N, x_{k}^N, m^N(\by^N)) i^N_k- \frac{1}{N^2}\sum_{k  = 1}^N \si^2: \pa_{\m \m}^2 \f_{1,1}(t^N, m^N(\by^N), x_{k}^N, x_{k}^N) i^N_k\nonumber \\
&&\dis = - [\dbL  \f + F] (t^N, m^N(\by^N)) - {1\over N}  \int_{\dbR^d} \si^2:\pa_{\m \m}^2 \f_{1,1}(t^N,m^N(\by^N),x,x) m^N(\by^N)(dx,1),\nonumber
\eea
where $\si (s, \tilde m, x, x) :=  \si (s, x, \tilde m)$. Note that, as $N\to\infty$,
\beaa
\int_{\dbR^d} \!\!\! \si^2:\pa_{\m \m}^2 \f_{1,1}(t^N,m^N(\by^N),x,x) m^N(\by^N)(dx,1)\to \int_{\dbR^d}  \!\!\! \si^2:\pa_{\m \m}^2 \f_{1,1}(0,m,x,x) m(dx,1),
\eeaa
since $(t^N, m^N(\by^N)) \to (0,m)$ and the mapping  $\dis(s, \tilde m) \mapsto \int_{\dbR^d}  \!\!\! \si^2:\pa_{\m \m}^2 \f_{1,1}(s,  \tilde m, x, x) \tilde m(dx,1)$ is continuous. Then \reff{subcv-phi}  and \reff{subcv-phi2} imply that
\beaa
 -(\dbL \f + F)(t^N, m^N(\by^N))  +   \cO(\frac 1 N)\le 0.
 \eeaa
Sending $N\to\infty$ and noticing that  $\dbL \f + F$ is continuous, we conclude the  viscosity subsolution property of $\ol U$ at $(0, m)$.

(ii) We now prove the  viscosity supersolution property. Fix $\f \in \ol \cA \ul U(0,m)$ 
with corresponding $\d_0$ and  s.t. $(0,m)$ is a strict maximizer of $\f - \ul U$ on $\cN_{\d_0}(0,m)$. With the same procedure and notations as in  (i) (in particular, we rewrite \eqref{estsup} with reversed inequalities and switch $\lim \sup$ and $\lim \inf$), we may find $(t_\d^N, m^N(\by_\d^N)) \underset{(\d,N) \to (0,\infty)}{\longrightarrow} (0,m)$ s.t. 
\beaa
 (\phi^N - U^N)(t_\d^N, \by_\d^N) = \max_{\th \in \cT_{t_\d^N, T}} \dbE_{t_\d^N, \by_\d^N}^{\dbP^0}\Big[(\phi^N - U^N)(\th \wedge \ch_\d^N, \bar \bY_{\th \wedge \ch_\d^N}) \Big].
 \eeaa
 Thus $\phi^N \in \ol \cA^N U^N(t_\d^N,  \by_\d^N)$ and the viscosity supersolution property provides
\bea\label{supercv1}
\left.\ba{lll}
\dis \min\Big\{\!\! -(\dbL \f + F)(t_\d^N, m^N(\by_\d^N)) + \cO\Big(\frac 1 N \Big), \ms\\
\dis \qq\q U^N(t_\d^N, m^N(\by_\d^N)) - \max_{\bi' < \bi_\d^N} U^N(t_\d^N, m^N(\bx_\d^N, \bi'))  \Big\} \ge 0. 
 \ea\right.
 \eea
As $(t_\d^N, m^N (\by_\d^N)) \underset{(\d, N) \to (0, \infty)}{\longrightarrow} (0,m)$,  sending $(\d, N) \longrightarrow (0,\infty)$ in \eqref{supercv1} provides the first part of the viscosity supersolution property:
\bea\label{supercv2}
-(\dbL  \f + F)(0, m)   \ge 0. \nonumber
\eea

It remains to show that $\ul U$ is nondecreasing for $\preceq$. We first assume that 
\bea\label{msmooth}
\mbox{$m(dx, 1)$ satisfies the required regularity in Lemma \ref{pure-mixed}.}
\eea
For any $A \in \cB(\dbR^d)$, by considering the corresponding cdf, we easily see that 
$$(m^N)^A(\by^N) := \frac 1 N \sum_{j=1}^N \d_{(x_k^N, i_k^N \1_{A}(x_k^N))} \underset{N \to \infty}{\longrightarrow} m^A.$$ 
By \eqref{supercv1}, we have $U^N(t^N, m^N(\by^N)) \ge U^N(t^N, (m^N)^A(\by^N))$ for all $N \ge 1$. Taking the $\lim \inf$ and using the fact that $U^N(t^N, m^N(\by^N)) \underset{N \to \infty}{\longrightarrow} \ul U(0,m)$, we deduce $\ul U(0,m) \ge \ul U(0,m^A)$. Given the arbitrariness of $A$ and by Lemma \ref{pure-mixed}, we deduce from the continuity of $\ul U$ that $\ul U(0,m) \ge \ul U(0,m')$ for all $m' \preceq m$, under the condition \eqref{msmooth}. 

Now for general $m$, there exist $m_n\in \cP_2(\bS)$  such that $\lim_{n\to\infty}W_2(m_n, m)=0$ and each $m_n$ satisfies \reff{msmooth}. Then $\ul U(0,m_n) \ge \ul U(0,m_n')$ for all $m'_n \preceq m_n$.  Recall \eqref{order} and let $m'\preceq m$ with transition function $p$. We first assume $p$ is continuous. Then, letting   $m_n'  \preceq m_n$ have the same $p$, we have $\lim_{n\to\infty}W_2(m_n', m')=0$. Thus, the continuity of $\ul U$ implies $\ul U(0,m) \ge \ul U(0,m')$.  Finally,  by \cite[Lemma 3.8]{TTZ2} we extend the inequality to all $m' \preceq m$, with possibly discontinuous $p$. 
\qed

 \section{Propagation of chaos for stopped  diffusions}
 \label{sec-proof-propagation}
 
 In this section we prove Theorems \ref{thm:propagation-optimal} and \ref{thm:approx}. We first use the latter to prove Theorem \ref{thm:propagation-optimal}. Recall \reff{weakoptstop}, \reff{Nparticles}, and that $\l$ denotes the canonical map on $\cP_2(\O)$. 
 
  \ms
 \no {\bf Proof of Theorem \ref{thm:propagation-optimal}.} (i) Let $\hat \n^N := \dbP^0 \circ (m^N(\bY^{\hat \btau^N}))^{-1}$. By Theorem  \ref{thm:approx} (i) and by otherwise considering a subsequence, we may assume $\hat\n^N$ converges weakly to some $\hat \n$ supported on $\cP(t,m)$.  For each $N\ge 1$, since $\hat \btau^N$ is $\e_N$-optimal for \eqref{Nparticles}, we have
 \bea
 \label{vNconv1}
 v^N(t, \by^N) \le J_N(t, \by^N, \hat\btau^N) + \e_N 
 = \dbE^{\dbP^0}\big[ J(t, m^N(\bY^{\hat \btau^N}))\big] + \e_N = \dbE^{\hat \n^N}[J(t,\l)] + \e_N.
 \eea
For any $R>0$, let $J_R(t, \l)$ be defined by truncating $F$ and $g$ with $R$ in \reff{weakoptstop}, and define $J_{N, R}$ similarly by truncation in \reff{Nparticles}. Then
\bea
\label{vNconv2}
\lim_{N\to\infty} \dbE^{\hat \n^N}[J_R(t,\l)] = \dbE^{\hat \n}[J_R(t,\l)].
\eea
By \reff{unif-estimate11}, we have the following uniform estimate:
\beaa
&&\dis \Big|\dbE^{\hat \n^N}[J_R(t,\l)]  - \dbE^{\hat \n^N}[J(t,\l)] \Big| = \Big|J_{N,R}(t, \by^N, \hat \btau^N)  - J_{N}(t, \by^N, \hat \btau^N)  \Big| \\
&&\dis \le C\dbE^{\dbP^0}\Big[{1\over N}\sum_{k=1}^N \int_t^T \big(1+ |X^{\hat\btau^N, k}_r| + \|\bX^{\hat\btau^N}_r\|_1\big) \1_{\{C(1+ |X^{\hat\btau^N, k}_r| + \|\bX^{\hat\btau^N}_r\|_1)\ge R\}} dr\\
&&\dis\qq  +  \|\bX^{\hat\btau^N}_T\|_1 \1_{\{C \|\bX^{\hat\btau^N}_T\|_1]\ge R\}}\Big]\\
&&\dis \le {C\over R}\dbE^{\dbP^0}\Big[{1\over N}\sum_{k=1}^N \int_t^T [1+ |X^{\hat\btau^N, k}_r|^2 + \|\bX^{\hat\btau^N}_r\|_1^2] dr  +  \|\bX^{\hat\btau^N}_T\|_1^2\Big]\\
&&\dis \le {C\over R}  \Big[{1\over N}\sum_{k=1}^N[1+|x^N_k|^2] + \|\bx^N\|_1^2\Big]  \le {C\over R} [1+ \|\bx^N\|_2^2].
\eeaa
Thus it is clear that
\beaa
\lim_{R\to\infty} \sup_{N\ge 1} \Big|\dbE^{\hat \n^N}[J_R(t,\l)]  - \dbE^{\hat \n^N}[J(t,\l)] \Big| = \lim_{R\to\infty}  \Big|\dbE^{\hat \n}[J_R(t,\l)]  - \dbE^{\hat \n}[J(t,\l)] \Big| =0.
\eeaa
Then by \reff{vNconv2} we have $\lim_{N\to\infty} \dbE^{\hat \n^N}[J(t,\l)] = \dbE^{\hat \n}[J(t,\l)]$. Send $N\to 0$ in \reff{vNconv1}, note that $\e_N\to 0$ and by Theorem \ref{thm-vfconv} $v^N(t, \by^N)  \underset{N \to \infty}{\longrightarrow} V(t,m)$, we obtain 
 $V(t,m) \le \dbE^{\hat \n}[J(t,\l)]$. 
Since $\hat \n$ is supported on $\cP(t,m)$, we have $V(t,m) = \sup_{\dbP \in \cP(t,m)} J(t,\dbP) \ge  J(t,\l)$, $\hat \n$-a.s.,  and therefore $V(t,m)=  J(t,\l)$, $\hat \n$-a.s. This proves that $\hat \n$ is supported on the set of optimal controls.
 
 (ii) By Theorem  \ref{thm:approx} (ii), there exists a sequence $\hat \btau^N$, $N \ge 1$, such that $m^N(\bY^{\hat \btau^N}) \overset{\cW_2}{\longrightarrow} \dbP^*$, $\dbP^0$-a.s. Then, by the continuity  of $F$ and $g$ and following similar arguments as in (i), 
 \bea
 \label{eNconv}
 \dbE^{\dbP^0}\big[J(t, m^N(\bY^{\hat \btau^N})) \big]\to J(t,\dbP^*)  = V(t,m).
 \eea
 Introduce 
$$ \e_N := V^N(t, m^N(\by^N)) - \dbE^{\dbP^0}\left[ J(t, m^N(\bY^{\hat \btau^N}))\right].
$$
By definition $\hat \btau^N$ is $\e_N$-optimal for \eqref{Nparticles}, and by Theorem \ref{thm-vfconv} and \reff{eNconv}, we get $\e_N\to 0$. 

(iii) Assume by contradiction that there exist $\e_N\to 0$ but $m^N(\bY^{\hat \btau^N}) $ does not converge to $\dbP^*$. Then there exists $c>0$ such that, by otherwise choosing a subsequence, $\dbP_0\big(\cW_2(m^N(\bY^{\hat \btau^N}), \dbP^*)>c\big) > c$. By Theorem \ref{thm:propagation-optimal} (i), again by otherwise choosing a subsequence, we may assume $\dbP^0 \circ \big(m^N(\bY^{\hat \btau^N})\big)^{-1}$ converges weakly to some $\nu\in \cP_2(\cP_2(\O))$ whose support is on the singleton $\{\dbP^*\}$, then we must have $m^N(\bY^{\hat \btau^N})\to \dbP^*$, in probability $\dbP^0$. This is a desired contradiction. 
\qed

 \ms
 \no {\bf Proof of Theorem \ref{thm:approx}.} Again for simplicity we assume $t=0$.
 
 (i) For each $N \ge 1$, denote $\n^N := \dbP^0 \circ (m^N(\bY^{\btau^N}))^{-1} \in \cP_2(\cP_2(\O))$.  
 
 \no{\bf Step 1.} We first show that $\{ \n^N \}_{N \ge 1}$ is uniformly integrable, that is, 
 \bea
 \label{ui}
 \lim_{R \to \infty} \sup_{N \ge 1} \dbE^{\n^N}\Big[ \|\l\|_2\1_{\{\|\l\|_2\ge R\}} \Big] = 0.
 \eea
 Indeed, for any $R > 0$, by \reff{unif-estimate11} we have
\beaa
&&\dis \dbE^{\n^N}\Big[  \|\l\|_2\1_{\{\|\l\|_2\ge R\}}  \Big] \le {1\over R} \dbE^{\n^N}\Big[ \|\l\|_2^2\Big] =\frac 1 R  \dbE^{\dbP^0}\Big[ \frac 1 N \sum_{k=1}^N \lvert \sup_{0\le t\le T}|Y^{\btau^N, k}_t|^2  \Big] \\
&&\dis \le \frac C R \Big[ 1+  \frac 1 N \sum_{k=1}^N [ |x^N_k|^2 + \|\bx^N\|_2^2] \Big] \le \frac C R \Big[ 1+   \|\by^N\|_2^2 \Big].
\eeaa
Since $\cW_2\big(m^N(\by^N), m\big) \underset{N \to \infty}{\longrightarrow} 0$, we have $\|\by^N\|_2\to \|m\|_2$, and thus $\{\|\by^N\|_2\}_{N\ge 1}\}$ is bounded. This implies \reff{ui} immediately.

\no{\bf Step 2.}  We next prove the tightness of $\{\nu^N\}_{N \ge 1}$.  By Lacker \cite[Corollary B.1]{LackerMartingale}, given \reff{ui}, it suffices to show that the mean measures $\{ \dbE^{\dbP^0}\big[m^N(\bY^{\btau^N}) \big]\}_{N \ge 1}\subset \cP_2(\O)$ is tight. Here, for a random measure  $\tilde m : \O^0 \to \cP_2(\O)$, the mean measure $\dbE^{\dbP^0}[\tilde m]$ is defined by
$$ \big\langle \dbE^{\dbP^0}[\tilde m], \f \big\rangle := \dbE^{\dbP^0} \big[ \langle \tilde m, \f \rangle \big] \q \mbox{for all $\f \in C_b^0(\O)$}. $$

To see this, we first remark that, the tightness of the joint measure is equivalent to the tightness of the marginals. For the first marginal, by \reff{unif-estimate12}  we have
{\small $$
\sup_{N \ge 1} \sup_{\t \in \cT_{0,T}} \Big\langle \dbE^{\dbP^0}\big[m^N(\bY^{\btau^N}) \big], \lvert X_{(\t+\d) \wedge T} - X_\t \rvert^2 \Big\rangle  =  \sup_{N \ge 1} \sup_{\t \in \cT_{0,T}} \frac 1 N \sum_{k=1}^N \dbE^{\dbP^0}\Big[ \big\lvert X_{(\t+\d) \wedge T}^{\btau^N,k} - X_\t^{\btau^N, k} \big\rvert^2 \Big] \le C\d.
$$ }
Then it follows from Aldous' criterion (see Billingsley \cite[Theorem 16.10]{Bil}) that the first marginal of $\{ \dbE^{\dbP^0}\big[m^N(\bY^{\btau^N}) \big]\}_{N \ge 1}$ is tight. Moreover, since the second marginals are measures on $\dbI^0([-1,T])$, which are in continuous bijection with measures on $[0,T]$, and thus are tight. Therefore, $\{ \dbE^{\dbP^0}\big[m^N(\bY^{\btau^N}) \big]\}_{N \ge 1}$ and hence $\{\nu^N\}_{N \ge 1}$ are tight.

\no{\bf Step 3.} We now show that $\n$ is supported on $\cP(0,m)$, or equivalently, that $\l \in \cP(0, m)$, $\n$-a.s. Since $\cW_2(m^N(\by^N), m) \underset{N \to \infty}{\longrightarrow} 0$, we have $\l_{0-} = m$, $\n$-a.s. Recall \reff{martingalepb}, we may consider equivalently the martingale problem. That is, for any $\psi \in C_b^2(\dbR^d)$, we want to show that
\bea
\label{mg}
M_s^{\psi} := \psi(X_s) - \int_0^s \cL_x \psi(r, X_r, \l_r)I_r  dr ~\mbox{is a $\l$-martingale, for $\n$-a.e. $\l$}.
\eea
For this purpose, fix $0\le s_1 < s_2\le T$ and $h_{s_1} \in C_b^0(\O)$ that is $\cF_{s_1}$-measurable. Note that  
 \beaa
&&\dis \dbE^{\nu^N}\Big[ \big|\la \l, h_{s_1}[M^\psi_{s_2}-M^{\psi}_{s_1}]\ra\big|\wedge 1\Big] \le \dbE^{\nu^N}\Big[ \big|\la \l, h_{s_1}[M^\psi_{s_2}-M^{\psi}_{s_1}]\ra\big|^2\Big]\\
&&\dis \hspace{-7mm} =\dbE^{\dbP^0}\Big[ \big|\la m^N(\bY^{\t^N}), h_{s_1}[M^\psi_{s_2}-M^{\psi}_{s_1}]\ra\big|^2\Big] \\
&&\dis \hspace{-7mm} =\dbE^{\dbP^0}\Big[ \Big|{1\over N}\sum_{k=1}^N h_{s_1}(Y^{\btau^N,k})\big[\psi(X^{\btau^N, k}_{s_2}) - \psi(X^{\btau^N, k}_{s_1})  - \!\! \int_{s_1}^{s_2} \!\!\! \cL_x \psi(r, X^{\btau^N,k}_r, m^N(\bY^{\btau^N}_r)) I^{\btau^N, k}_r  dr\big]\Big|^2\Big].
\eeaa
By \reff{Ndynamics} and \reff{symmetric}, it follows from the It\^{o} formula that
\beaa
&&\dis \dbE^{\nu^N}\Big[ \big|\la \l, h_{s_1}[M^\psi_{s_2}-M^{\psi}_{s_1}]\ra\big|\wedge 1\Big] \\
&&\dis\le \dbE^{\dbP^0}\Big[ \Big|{1\over N}\sum_{k=1}^N h_{s_1}(Y^{\btau^N,k}) \int_{s_1}^{s_2} \pa_x \psi( X^{\btau^N,k}_r) \si (r,  X^{\btau^N,k}_r, m^N(\bY^{\btau^N}_r)) I^{\btau^N,k}_r  dW_r\Big|^2\Big]\\
&&\dis = {1\over N^2}\sum_{k=1}^N \dbE^{\dbP^0}\Big[  \big|h_{s_1}(Y^{\btau^N,k})\big|^2 \int_{s_1}^{s_2} \big|\pa_x \psi( X^{\btau^N,k}_r) \si (r,  X^{\btau^N,k}_r, m^N(\bY^{\btau^N}_r))I^{\btau^N,k}_r \big|^2 dr\Big]\\
&&\dis \le {C\over N^2}\sum_{k=1}^N \dbE^{\dbP^0}\Big[   \int_{s_1}^{s_2} \big| \si (r,  X^{\btau^N,k}_r, m^N(\bY^{\btau^N}_r))\big|^2 dr\Big] \le {C\over N},
\eeaa
where the last inequality thanks to \reff{unif-estimate11}.  Note that $\nu^N$ converges to $\nu$ weakly, and the mapping $\l \in \cP_2(\O) \longmapsto  \big|\la \l, h_{s_1}[M^\psi_{s_2}-M^{\psi}_{s_1}]\ra\big|\wedge 1$ is bounded and continuous. Then
\beaa
\dbE^{\nu}\Big[ \big|\la \l, h_{s_1}[M^\psi_{s_2}-M^{\psi}_{s_1}]\ra\big|\wedge 1\Big] =0.
\eeaa
This implies that, for any desired $\psi$, $h_{s_1}$ and $s_1<s_2$,
\bea
\label{mg2}
\la \l, h_{s_1}[M^\psi_{s_2}-M^{\psi}_{s_1}]\ra = 0, \q\mbox{for $\n$-a.e. $\l$}.
\eea
Since $M^\psi$ is continuous in $s$, and the spaces of the above $\psi$ and $h_{s_1}$ are separable, then \reff{mg2} implies \reff{mg}. Therefore, $\l\in \cP(0, m)$, $\nu$-a.s.

(ii) Without loss of generality, assume $(\O^0, \cF_T^0, \dbF^0, \dbP^0)$ is sufficiently large so that there exist a sequence of i.i.d.\ processes $\{(I^{0, k}, W^k)\}_{k\ge 1}$ which are independent of $\cF^0_{0-}$ (under $\dbP^0$) and $\dbP^0_{(I^{0,k}, W^k)} = \dbP_{(I, W^\dbP)|I_{0-} = 1}$. In particular, $W^k$ is a $\dbP^0$-Brownian motion and $I^{0,k}_{0-}=1$. Now for given $\by^N$, define
\bea
\label{tauN}
\dis I^{N, k}_t := I^{0, k}_t \1_{\{i_k = 1\}};\q \t^N_k := \inf\{t\ge 0: I^{N, k}_t = 0\},\q \btau^N:=(\t^N_1,\cds, \t^N_N).
\eea
We shall prove in four steps that the above $\btau^N$ satisfies the convergence requirement. Since these stopping times are adapted to the filtration $\dbF^0$ which might be larger than the Brownian filtration, we finally show in Step 5 that they can be approximated by pure stopping strategies, and that we can derive from this approximation the desired result. 

\no{\bf Step 1.} We first assume
\bea
\label{3integrable}
\int_\bS |x|^3 m(dy) <\infty,\q\mbox{and}\q \lim_{N\to\infty} \cW_3(m^N(\by^N), m) =0.
\eea
Denote $m_t := \dbP_{Y_t}$ and introduce
\bea
\label{tildeXNk}
\tilde X^{N, k}_t = x^N_k + \int_0^t b(r,\tilde X^{N, k}_r, m_r) I^{N, k}_r dr + \int_0^t \si(r, \tilde X^{N, k}_r, m_r) I^{N, k}_r dW^k_r.
\eea
Note that $(W^k, I^{N, k}, \tilde X^{N, k})$, $k \in [N]$, are independent. In this step we show that
\bea
\label{tildemconv}
\lim_{N\to\infty} \sup_{0\le t\le T} \dbE^{\dbP^0}\Big[\cW_1^2(m^N(\tilde \bY^N_t), m_t)\Big]=0, \q \mbox{for all $t$, where}\q \tilde \bY^{N, k} := (\tilde X^{N, k}, I^{N, k}).
\eea
To see this,  let $E^N_k\in \cF^0_{0-}$, $k \in [N]$ be a partition of $\O^0$ such that $\dbP^0(E^N_k)={1\over N}$. Denote
\bea
\label{hatXN}
(\hat W^N, \hat X^N_t, \hat I^N_t) := \sum_{k=1}^N   (W^k, \tilde X^{N, k}_t, I^{N, k}_t) \1_{E^N_k},\q \hat Y^N:=(\hat X^N_t, \hat I^N_t).
\eea
Then it follows from the arguments in \cite[Lemma 8.4]{MZ} that 
\bea
\label{hatXNest1}
\dbE^{\dbP^0}\Big[\cW_1^2(m^N(\tilde\bY^N_t), \dbP^0_{\hat Y^N_t})\Big] \le {C\over N^{2\over (d+1)\vee 3}} \Big(\dbE^{\dbP^0}[|\hat Y^N_t|^3]\Big)^{2\over 3}.
\eea
We shall point out that \cite[Lemma 8.4]{MZ} assumes Brownian filtration, which is not the case here. However, we emphasize that this assumption is due to the setting in \cite{MZ} and is never used there. In fact, the arguments provided in the proof of \cite[Lemma 8.4]{MZ} used only the independence between  $\{(\tilde X^{N, k}, I^{N, k})\}_k$  and $\{E^N_k\}_k$, which holds true here. Then \cite{MZ} refers to  \cite[Lemma 5, Lemma 6, and Theorem 1]{FournierGuillin}, which do not require the Brownian filtration.
Moreover, by standard SDE estimate and noting that $\cW_3(m^N(\by^N), m)\to 0$, we may estimate $\dbE^{\dbP^0}[|\hat Y^N_t|^3]$ and thus derive from \reff{hatXNest1} that
\bea
\label{hatXNest2}
\dbE^{\dbP^0}\Big[\cW_1^2(m^N(\tilde\bY^N_t), \dbP^0_{\hat Y^N_t})\Big]  \le {C(1+ \|\bx^N\|_3^2)\over N^{2\over (d+1)\vee 3}} \le{C_{m} \over N^{2\over (d+1)\vee 3}}, 
\eea
for some constant $C_{m}$, which may depend on $m$, but not on $N$. 

Next, for any bounded function $\f\in C^0(\O)$. Denote 
\beaa
v(i) := \dbE^\dbP[ \f(W^\dbP, I) | I_{0-} = i],\q i=0,1.
\eeaa
Note that $(W^k, \tilde X^{N, k}_t, I^{N, k}_t)$ is independent of $\cF^0_{0-}$. Then, by \reff{tauN},
\bea
\label{weakconv}
&&\dis \dbE^{\dbP^0}[ \f(\hat W^N, \hat I^N)] = {1\over N} \sum_{k=1}^N\dbE^{\dbP^0}[ \f(W^k, I^{N, k})] = {1\over N} \sum_{k=1}^N\Big[v(1) \1_{\{i_k=1\}} + v(0) \1_{\{i_k=0\}} \Big]\nonumber\\
&&\dis \longrightarrow v(1) \dbP(I_{0-} =1 ) + v(0) \dbP(I_{0-}=0) = \dbE^\dbP[\f(W^\dbP, I)],
\eea
as $N\to\infty$. That is, $\dbP^0\circ (\hat W^N, \hat I^N)^{-1} \to \dbP\circ (W^\dbP, I)^{-1}$ weakly. Since $I$ is bounded, by \cite[Theorem 5.5]{CarDel} we have
\beaa
\lim_{N\to\infty} \cW_2\Big(\dbP^0\circ (\hat W^N, \hat I^N)^{-1},~ \dbP\circ (W^\dbP, I)^{-1}\Big) =0.
\eeaa
Moreover, since $\lim_{N\to\infty} \cW_2\big(m^N(\by^N), m)\big) =0$, then
\bea
\label{weakconv2}
\lim_{N\to\infty} \cW_2\Big(\dbP^0\circ (\hat X^N_0, \hat W^N, \hat I^N)^{-1},~ \dbP\circ (X_0, W^\dbP, I)^{-1}\Big) =0. \nonumber
\eea
One may easily verify that
\beaa
\hat X^N_t = \hat X^N_0 + \int_0^t b(r,\hat X^{N}_r, m_r) \hat I^{N}_r dr + \int_0^t \si(r, \hat X^{N}_r, m_r) \hat I^{N}_r d\hat W^N_r.
\eeaa
Then, by comparing this with \reff{asympt}, 
$\lim_{N\to\infty} \sup_{0\le t\le T} \cW_2\big(\dbP^0_{\hat Y^N_t}, m_t\big)  = 0.
$
This, together with \reff{hatXNest2}, implies \reff{tildemconv} immediately. 

\no{\bf Step 2.} In this step we show that, again under the additional condition \reff{3integrable},
\bea
\label{taumconv}
\lim_{N\to\infty} \sup_{0\le t\le T} \dbE^{\dbP^0}\Big[\cW_1^2(m^N(\bY^{\btau^N}_t), m_t)\Big]=0.
\eea
Recall \reff{tauN} and \reff{symmetric}, and compare \reff{Ndynamics} and \reff{tildeXNk}. By the Lipschitz continuity of $b, \si$, especially the $\cW_1$-Lipschitz continuity in $m$, it follows from standard SDE estimates that
\beaa
 \dbE^{\dbP^0}\Big[|X^{\btau^N, k}_t-\tilde X^{N, k}_t |^2\Big] \le C\dbE^{\dbP^0}\Big[\int_0^t \cW_1^2( m^N(\bY^{\t^N}_s), m_s) ds\Big].
\eeaa
Then
\beaa
 &&\dis \dbE^{\dbP^0}\Big[\cW_1^2(m^N(\bY^{\btau^N}_t), m_t)\Big] \le C \dbE^{\dbP^0}\Big[\cW_1^2(m^N(\bY^{\btau^N}_t), m^N(\tilde \bY^N_t)) + \cW_1^2(m^N(\tilde \bY^N_t), m_t))\Big]\\
 &&\dis \le C \dbE^{\dbP^0}\Big[ {1\over N} \sum_{k=1}^N |X^{\btau^N,k}_t-\tilde X^{N, k}_t |^2 + \cW_1^2(m^N(\tilde \bY^N_t), m_t))\Big]\\
 &&\dis \le C\dbE^{\dbP^0}\Big[ \int_0^t \cW_1^2( m^N(\bY^{\t^N}_s), m_s) ds+ \cW_1^2(m^N(\tilde \bY^N_t), m_t))\Big]
 \eeaa
 Apply the Gronwall inequality, we obtain
 \beaa
  \dbE^{\dbP^0}\Big[\cW_1^2(m^N(\bY^{\btau^N}_t), m_t)\Big] \le  C\dbE^{\dbP^0}\Big[\cW_1^2(m^N(\tilde \bY^N_t), m_t))\Big].
  \eeaa
  This, together with \reff{tildemconv}, implies \reff{taumconv} immediately.
 
 \no{\bf Step 3.} We now prove the result without assuming \reff{3integrable}. For any $R>0$, let $\phi_R(x)$ denote the truncation function such that $x$ is truncated by $R$. Denote $m^R:= m\circ (\phi_R(X_0), I_{0-})^{-1}$, $\by^{N,R} := (\phi_R(x^N_k), i^N_k)_{1\le k\le N}$.  Then it is clear that
 \bea
 \label{MRconv1}
&\dis \lim_{R\to\infty} \cW_2(m^R, m) =0,\q \lim_{R\to\infty} \sup_{N\ge 1} {1\over N}\sum_{k=1}^N \big|\phi_R(x^N_k)- x^N_k\big|^2=0,\\
 \label{MRconv2}
 &\dis \dis \lim_{N\to\infty} \cW_3(m^N(\by^{N,R}), m^R)=0,\q \mbox{for all} \ R>0.
\eea
Introduce $\dbP^R$ and $\bY^{R,\btau}$ in an obvious way. By \reff{taumconv} and \reff{MRconv2} we have 
\bea
\label{MRconv3}
\lim_{N\to\infty} \sup_{0\le t\le T}\dbE^{\dbP^0}\Big[\cW_1^2(m^N(\bY^{R,\btau^N}_t), \dbP^R_{Y_t})\Big]=0.
\eea
Moreover,  by standard SDE estimates one can easily see that
\beaa
&&\dis \cW_2( \dbP^R_{Y_t}, m_t) \le C\cW_2(m^R, m),\\
&&\dis \dbE^{\dbP^0} \Big[\cW_2^2\big(m^N(\bY^{R, \btau^N}_t), m^N(\bY^{\btau^N}_t)\big)\Big] \le {1\over N}\sum_{k=1}^N \dbE^{\dbP^0} \Big[\big|X^{R, \btau^N,k}_t-X^{\btau^N,k}_t\big|^2\Big] 
\\&&
\dis \hspace{60mm} 
\le {C\over N}\sum_{k=1}^N \big|\phi_R(x^N_k)- x^N_k\big|^2.
\eeaa
Then it follows from \reff{MRconv1}  that, for any $R>0$,
\beaa
\dis\underset{N\to\infty}{\lim \sup}\sup_{0\le t\le T}\dbE^{\dbP^0} \Big[\cW_1^2\big( m^N(\bY^{\btau^N}_t), m_t\big)\Big] \!\! &\le& \!\!  \underset{N\to\infty}{\lim \sup} \ C\sup_{0\le t\le T}\dbE^{\dbP^0} \Big[\cW_1^2\big( m^N(\bY^{\btau^N}_t), m^N(\bY^{R, \btau^N}_t)) \\
&&\dis\q + \cW_1^2(m^N(\bY^{R,\btau^N}_t), \dbP^R_{Y_t}) +\cW^2_2( \dbP^R_{Y_t}, m_t)\big)\Big]\\
 \!\!& =&  \!\! C\sup_{0\le t\le T}\dbE^{\dbP^0} \Big[ \cW_1^2(m^N(\bY^{R,\btau^N}_t), \dbP^R_{Y_t}) \Big].
\eeaa
This, together with \reff{MRconv3}, implies \reff{taumconv}, without assuming \reff{3integrable}.
 
\no{\bf Step 4.} Finally we prove the convergence of the distribution of the processes. First, compare \reff{Ndynamics} and \reff{tildeXNk} again, by \reff{taumconv} we obtain immediately that
\bea
\label{jointconv1}
\lim_{N\to\infty}\dbE^{\dbP^0}\Big[\sup_{0\le t\le T} |X^{\btau^N,k}_t-\tilde X^{N, k}_t |^2\Big] =0.
\eea
Next,  for any bounded function $\f\in C^0(\O)$, similarly to \reff{weakconv}, we can show  
\beaa
\lim_{N\to\infty} \dbE^{\dbP^0}\Big[\dbE^{m^N(\tilde \bY^N)}[\f(Y)] \Big] = \lim_{N\to\infty} \dbE^{\dbP^0}\Big[{1\over N} \sum_{k=1}^N [\f(\tilde Y^{N, k})] \Big] = \dbE^\dbP[\f(Y)].
\eeaa
This, together with \reff{jointconv1}, implies
\bea
\label{jointconv2}
\lim_{N\to\infty} \dbE^{\dbP^0}\Big[\dbE^{m^N(\bY^{\btau^N})}[\f(Y)] \Big]  = \dbE^\dbP[\f(Y)].
\eea
On the other hand, by (i) of this theorem, we know $\big\{\dbP^0 \circ \big(m^N(\bY^{\btau^N})\big)^{-1} \big\}_{N \ge 1}$ is tight. By \reff{jointconv2} clearly $\dbP$ is its unique accumulation point. Then, similarly to the arguments in Theorem \ref{thm:propagation-optimal} (iii), we obtain $m^N(\bY^{\btau^N})  \overset{\cW_2}{\longrightarrow} \dbP$, $\dbP^0$-a.s.

\no {\bf Step 5:} It remains to approximate the sequence $\{\btau^N\}_{\{N \ge 1\}}$ with a sequence of $\tilde \btau^N \in \mathbfcal{T}_{t,T}^N$, that is, of stopping times adapted to the Brownian filtration. Indeed, the stopping times constructed in the previous steps are adapted to a larger filtration. Fixing $N \ge 1$, we have by Carmona, Delarue \& Lacker \cite[Theorem 6.4]{CDL} the existence of a sequence $\{\tilde \btau^{N, k}\}_{\{k \ge 1\}}$ in $\mathbfcal{T}_{t,T}^N$ such that $(\tilde \btau^{N,k}, \boldsymbol{W}^N$ converges in distribution to $(\btau^N, \boldsymbol{W}^N)$ as $k \to \infty$, with $\boldsymbol{W}^N := (W^1, \dots, W^N)$. This implies that $\bY^{N, \tilde \btau^{N,k}}$ converges in distribution to $\bY^{N, \btau^{N}}$, and therefore that $m^N\big(\bY^{N, \tilde \btau^{N,k}})$ converges in distribution to $m^N\big(\bY^{N, \btau^{N}})$. 

It is also obvious from the previous steps that $m^N\big(\bY^{N, \btau^{N}})$ converges in distribution to $\dbP$ as $N \to \infty$. Then there exists a subsequence such that $m^N\big(\bY^{N, \tilde \btau^{N,k_N}})$ converges in distribution to $\dbP$ as $N \to \infty$. Note that $\dbP$ is a constant in the space of measure-valued random variables (endowed with the distance $\cW_2$); therefore the previous convergence also holds in probability. We may then extract a subsequence (still denoted the same) converging to $\dbP$, $\dbP^0$-a.s., that is,
$$ \cW_2\big(m^N\big(\bY^{N, \tilde \btau^{N,k_N}}), \dbP\big) \underset{N \to \infty}{\longrightarrow} 0, \ \mbox{$\dbP^0$-a.s.} $$
This concludes the proof.

\qed

\section{Further study on the multiple optimal stopping problem}
\label{sect-Ndimpb}

\subsection{The viscosity solution property}
In this subsection we prove Theorem \ref{thm:Nexistence}. Recall the notations in Subsection \ref{sect-general} and \ref{sec:viscosityN}, in particular the unstopped process $\bar \bY$.  
As standard, the viscosity property relies on the following dynamic programming principle. 

\begin{prop}\label{DPPKQR}
Under Assumption \ref{assum-N}, for all $(t,\by) \in \L^N_0$, we have
\bea\label{Ntostd}
v^N(t,\by) = \sup_{\th \in \cT_{t,T}^N} \dbE_{t,\by}\Big[\frac 1 N \sum_{k=1}^N \int_t^\th f^k(r, \bar \bX_r)\bar I_r^kdr + \max_{\bi' < \bi} v^N(\th, \bar \bX_\th, \bi') \Big].
\eea
In particular, by setting $\th = t$, we see that $v^N$ is non-decreasing in $\bi$:
\bea
\label{vNmon}
v^N(t,\by)  \ge v^N(t, \bx, \bi'),\q\mbox{for all}\q \bi'<\bi.
\eea

\end{prop}
\proof First, by Kobylanski, Quenez \& Rouy-Mironescu \cite[Theorem 3.1]{KQR}\footnote{Although Kobylanski, Quenez \& Rouy-Mironescu assume nonnegative reward processes throughout their paper, the validity of their result in our context is easy to check.}, we have
\beaa
v^N(t,\by) &=& \sup_{\th \in \cT_{t,T}} \dbE_{t,\by}\Big[ \max_{l \in [N]} \sup_{\btau \in \cT_{\th, T}^{N-1}} \dbE_{\th, \bY_\th^{-l}}\Big[\frac 1 N \sum_{k=1}^N \int_t^T f^k(r, \bX_r^{\btau \otimes_l \th})\bar I_r^kdr + g\big(\bX_T^{\btau \otimes_l \th}\big) \Big] \Big] \\
&=& \sup_{\th \in \cT_{t,T}} \dbE_{t,\by}\Big[ \frac 1 N \sum_{k=1}^N \int_t^\th f^k(r, \bar \bX_r)\bar I_r^kdr \\
&& \q + \max_{l \in [N]} \sup_{\btau \in \cT_{\th, T}^{N-1}} \dbE_{\th, \bY_\th^{-l}}\Big[\frac 1 N \sum_{k=1}^N \int_\th^T f^k(r, \bX_r^{\btau \otimes_l \th})\bar I_r^kdr + g\big(\bX_T^{\btau \otimes_l \th}\big) \Big] \Big],
\eeaa
where 
\bea
\label{bi-l}
\btau \otimes_l \th := (\t_1, \dots, \t_{l-1}, \th, \t_{l+1}, \dots, \t_N),\q \bi^{-l} := (i_1, \dots, i_{l-1}, 0, i_{l+1}, \dots, i_N),
\eea 
and $\bY^{-l} := (\bar \bX, \bi^{-l})$.
By the Markov property, we have 
$$ 
\sup_{\btau \in \cT_{\th, T}^{N-1}} \dbE_{\th, \bY_\th^{-l}}\Big[\frac 1 N \sum_{k=1}^N \int_\th^T f^k(r, \bX_r^{\btau \otimes_l \th})\bar I_r^kdr + g\big(\bX_T^{\btau \otimes_l \th}\big) \Big] = v^N(\th, \bar \bX_\th, \bi^{-l}). 
$$
Finally, by induction $v^N$ is nondecreasing in $\bi$ for the partial order on $\{0,1\}^N$, and thus
$$ 
\max_{l \in [N]}  v^N(\th, \bar \bX_\th, \bi^{-l}) = \max_{\bi' < \bi}  v^N(\th, \bar \bX_\th, \bi').   
$$

\vspace{-7mm}
\qed

\no {\bf Proof of Theorem \ref{thm:Nexistence}.}   We proceed by induction on the cardinality $|A(\bi)|$ of the set
\bea
\label{Ai}
A(\bi):=\{k\in [N]: i_k = 1\}.
\eea
As $(t, \bx)\longmapsto g_N(t, \bx)$ is continuous, it suffices to prove that whenever $v^N(t, \bx, \bi)$ is continuous in $(t, \bx)$ for $|A(\bi)| \le n$, we may conclude for $|A(\bi)| =n+1$ that $v^N(\cdot, \bi)$ is the unique continuous viscosity solution of \eqref{cascade}. Note that $\bar I^k_r \equiv 1$ for $k\in A(\bi)$ and $\bar I^k_r \equiv 0$ for $k\notin A(\bi)$. To see this, we rewrite \reff{Ntostd}, emphasizing again that the process $\bar \bX$ is unstopped, thus reducing \reff{Ntostd2} to a standard optimal stopping problem: 
\begin{equation}
\label{Ntostd2} 
v^N(t,\bx, \bi) = \!\sup_{\th \in \cT_{t,T}^N} \dbE_{t,\by}\Big[\frac 1 N \!\!\sum_{k\in A(\bi)} \!\int_t^\th \!\!f^k(r, \bar \bX_r)dr + \hat v^N\!(\th, \bar \bX_\th, \bi) \Big],
~~
\hat v^N(t, \bx, \bi) := \sup_{\bi'<\bi}  v^N(t, \bx, \bi').
\end{equation}
Notice that $\hat v^N(t, \bx, \bi)$ is continuous in $(t, \bx)$, as $|A(\bi')| \le n$. Then it follows from  Ekren \cite{Ekren} that $v^N(\cd, \cd, \bi)$ is a continuous viscosity solution of \reff{cascade} in the sense of  Definition \ref{defn:viscosityN}.\footnote{Here again, we notice that our conditions are slightly different from Ekren and that the viscosity solution property is easily checked to be valid in our context.} Moreover, note that a viscosity  sub/super solution in the sense of Definition \ref{defn:viscosityN} implies it is a viscosity  sub/super solution in the classical Crandall-Lions sense.  Then the comparison principle for viscosity semi-solutions in Crandall-Lions sense implies that in our sense, see Remark \ref{rem-viscosityN}, and hence $v^N$ is the unique viscosity solution to \reff{cascade}.
\qed

\subsection{The symmetric case}
We now prove Lemma \ref{lem:unif-estimate} and Theorem \ref{thm:Nexistence-symmetric} in the symmetric case \reff{symmetric}.

\ms
\no{\bf Proof of Lemma \ref{lem:unif-estimate}.}  (i) Recall \reff{Ndynamics}. For $k \in [N]$, by standard SDE estimates we have
\bea
\label{unif-est1}
\dbE^{\dbP^0}[|X^{\btau, k}_s|^2]  \le \dbE^{\dbP^0}[\sup_{t\le r\le s}|X^{\btau, k}_r|^2] \le C\dbE^{\dbP^0}\Big[ |x_k|^2 + \int_t^s \big(1+ \|\bX^\t_r\|_2^2\big) dr\Big].
\eea
Then
\beaa
\dbE^{\dbP^0}[\|\bX^\btau_s\|_2^2] \le C\dbE^{\dbP^0}\Big[ \|\bx\|_2^2 +1+ \int_t^s  \|\bX^{\btau}_r\|_2^2 dr\Big].
\eeaa
Applying Gronwall inequality, we have $ \dbE^{\dbP^0}\big[\sup_{t\le s\le T}\|\bX^\btau_s\|_2 \le C(1+ \|\bx\|_2)\big]$. Plug this into the second inequality of \reff{unif-est1}, we obtain  \reff{unif-estimate11} immediately.
 
 Moreover, for any $\th_1, \th_2\in \cT^N_{[t, T]}$ with $\th_1\le \th_2$, by standard SDE estimates we have
\beaa
&&\dis \dbE^{\dbP^0}[\sup_{\th_1\le s\le \th_2}|X^{\btau, k}_s-X^{\btau, k}_{\th_1}|^2] \le C\dbE^{\dbP^0}\Big[ \int_{\th_1}^{\th_2}\big( 1+  |X^{\btau, k}_r|^2 + \|\bX^\btau_r\|^2\big) dr\Big] \\
&&\le C\|\th_2-\th_1\|_\infty \dbE^{\dbP^0}\Big[  1+ \sup_{t\le r\le T} |X^{\btau, k}_r|^2 + {1\over N}\sum_{l=1}^N\sup_{t\le r\le T} |X^{\btau, l}_r|^2\Big].
\eeaa
Plug \reff{unif-estimate11} into it, we obtain \reff{unif-estimate12}.

(ii) We now assume Assumption \ref{assum-comparison} (i) also holds true.  Let $\d>0$ be a small number which will be specified later, and consider a partition $t=t_0<\cds<t_n=T$ such that $\D t_j \le \d$ for all $j$. For each $j=1,\cds, n$ and $k=1, \cds, N$, by standard SDE estimates we have, for a generic constant $C$ which is independent of $N$ and $\d$,
\beaa
\dbE^{\dbP^0}\Big[\sup_{t_j \le s\le t_{j+1}} |X^{\btau, k}_s|\Big]  &\le& C\dbE^{\dbP^0}\Big[|X^{\btau,k}_{t_j}| + \int_{t_j}^{t_{j+1}}\big(1+|X^{\btau, k}_s| + \|\bX^\btau_s\|_1\big) ds \\
&&+ \big(\int_{t_j}^{t_{j+1}}\big(1+|X^{\btau, k}_s| + \|\bX^\btau_s\|_1\big)^2ds\big)^{1\over 2}\Big]\\
&\le& C\dbE^{\dbP^0}\Big[|X^{\btau,k}_{t_j}| + 1+ \sqrt{\d} \sup_{t_j\le s\le t_{j+1}} [|X^{\btau, k}_s| + \|\bX^\btau_s\|_1] \Big]
\eeaa
Set $\d := {1\over 9 C^2}$ for the above $C$,  we have
\bea
\label{unif-est2}
{2\over 3}\dbE^{\dbP^0}\Big[\sup_{t_j \le s\le t_{j+1}} |X^{\btau, k}_s|\Big] \le C\dbE^{\dbP^0}\Big[|X^{\btau,k}_{t_j}| + 1\Big]+ {1\over 3}\dbE^{\dbP^0}\Big[ \sup_{t_j\le s\le t_{j+1}}  \|\bX^\btau_s\|_1\Big].
\eea
This implies that
\beaa
{2\over 3}\dbE^{\dbP^0}\Big[{1\over N}\sum_{k=1}^N\sup_{t_j \le s\le t_{j+1}} |X^{\btau, k}_s|\Big]  \le C\dbE^{\dbP^0}\Big[\|\bX^{\btau}_{t_j}\|_1 + 1\Big]+ {1\over 3}\dbE^{\dbP^0}\Big[ \sup_{t_j\le s\le t_{j+1}}  \|\bX^\btau_s\|_1\Big],
\eeaa
and thus
\beaa
\dbE^{\dbP^0}\Big[\sup_{t_j \le s\le t_{j+1}} \| \bX^{\btau}_s\|_1\Big]\le \dbE^{\dbP^0}\Big[{1\over N}\sum_{k=1}^N\sup_{t_j \le s\le t_{j+1}} |X^{\btau, k}_s|\Big]  \le C\dbE^{\dbP^0}\Big[\|\bX^{\btau}_{t_j}\|_1 + 1\Big].
\eeaa
By induction on $j=0,\cds, n-1$, one can easily obtain 
\beaa
\dbE^{\dbP^0}\Big[\sup_{t \le s\le T} \| \bX^{\btau}_s\|_1\Big]  \le C(1+\|\bx\|_1).
\eeaa
Plugging this into \reff{unif-est2} and by induction on $j=0,\cds, n-1$ again, we prove  \reff{unif-estimate21}. Moreover, \reff{unif-estimate22} follows similar arguments as for \reff{unif-estimate12}.
\qed

\ms

\no{\bf Proof of Theorem \ref{thm:Nexistence-symmetric}.}  It is obvious that Assumptions \ref{assum-bsig}, \ref{assum-comparison} and \reff{symmetric} imply Assumption \ref{assum-N}. In particular, the uniform continuity of $g^N$ comes from the fact that $g$ is uniformly continuous on all bounded parts of $\cP_2(\dbR^d)$, which are $\cW_1$-compact. Then (i) follows directly from Theorem \ref{thm:Nexistence}.

To see (ii), assume $m_N = m^N(\by)$ and recall \reff{Nparticles} and \reff{VNsymmetric}. By the linear growth of $f$ and $g$ and by \reff{unif-estimate11}, for any $\btau$ we have
\beaa
\dis \big|J_N(t, \by, \btau)\big|&\le& \dbE^{\dbP^0}\Big[  {C\over N} \sum_{k=1}^N \int_t^T \big(1+ |X^{\btau, k}_r| + \|\bX^\btau_r\|_1\big) dr + 1+ \|\bX^\btau_T\|_1\Big]\\
&\le& C+ {C\over N} \sum_{k=1}^N \dbE^{\dbP^0}\big[\sup_{t\le  r\le T}\|\bX^{\btau}_r\|_1\big] \le  C(1+\|\bx\|_1) \le C(1+\|m_N\|_1).
\eeaa
 Then by the arbitrariness of $\btau$ we prove \reff{unif-estimate31}.

It remains to prove (iii). We proceed in three steps.

\no{\bf Step 1.} By Birkhoff's theorem, the Kantorovitch solution of the optimal transport problem in the finite discrete case is a Monge solution, i.e.\ the optimal coupling measure is concentrated on a graph, see e.g.\ Villani \cite{Villani} p5. Then, for any $t\in [0, T]$, $\by, \tilde \by \in \bS^N$:
$$ \cW_2^2(m^N(\by), m^N(\tilde \by)) = \min_{\pi \in \mathfrak{S}_N} \frac 1 N \sum_{k=1}^N \lvert y_k - \tilde y_{\pi(k)} \rvert^2, $$
where $\mathfrak{S}_N$ is the set of permutations on $[N]$. By possibly reordering the $y_k$'s, we may then assume without loss of generality that
\bea
\label{Birkhoff}
\cW_2(m^N(\by), m^N(\tilde \by)) = \|\by-\tilde\by\|_2.
\eea
Given $\btau$, let $\bY^\btau$ and $\tilde \bY^\btau$ be defined by \reff{Ndynamics} with initial condition $\by$ and $\tilde \by$, respectively, and denote $\Delta \bY^\btau := \bY^\btau-\tilde \bY^\btau$, $\Delta\by :=\by-\tilde\by$. Recall that $\btau$ is adapted to the Brownian filtration, namely in the spirit of open loop controls, then 
\bea
\label{DItau}
\Delta I^{\btau, k}_s = \Delta i_k \1_{\{\t_k>s\}}\q\mbox{and thus}\q \|\D \bI^\btau_t\|_2 \le \|\D \bi\|_2,\q \forall t\in [0, T].
\eea
In this step we shall prove further that
\bea
\label{unif-est4}
 \dbE^{\dbP^0}\Big[{1\over N}\sum_{k=1}^N\sup_{t \le s\le T} |\D X^{\btau, k}_s|\Big] 
 \le C(1+\|\bx\|_2)\|\D\by\|_2.
\eea
Indeed, when $i_k=1$ and $\tilde i_k=0$, we have $\tilde \bX^{\btau, k}_s \equiv  \tilde x_k$, and thus, by  \reff{unif-estimate22},
\bea
\label{unif-est3}
\dbE^{\dbP^0}\Big[\sup_{t\le s\le T} |\D X^{\btau, k}_s|\Big] &\le& C\dbE^{\dbP^0}\Big[|\D x_k| + \sup_{t\le s\le T}|X^{\btau, k}_s-x_k|\Big] \nonumber\\
&\le& C(1 + |x_k| +\|\bx\|_1 +|\D x_k|).
\eea
Similarly, when $i_k=0$ and $\tilde i_k=1$, the above still holds true. Moreover,  when $\D i_k=0$, we have $\D I^{\btau, k}_s \equiv 0$. Let $\d$ and $t=t_0<\cds<t_n=T$ be as in the proof of Lemma \ref{lem:unif-estimate} (ii). Then, for each $j=0,\cds, n-1$, similar to \reff{unif-est2} we have
\beaa
{2\over 3}\dbE^{\dbP^0}\Big[\sup_{t_j \le s\le t_{j+1}} |\D X^{\btau, k}_s|\Big] \le {1\over 3}\dbE^{\dbP^0}\Big[ \sup_{t_j\le s\le t_{j+1}}  \|\D \bX^\btau_s\|_1\Big] +  C\dbE^{\dbP^0}\Big[|\D X^{\btau,k}_{t_j}| \Big].
\eeaa
This, together with \reff{unif-est3} which is for all $k$ such that $\D i_k\neq 0$, implies that
\beaa
 {2\over 3}\dbE^{\dbP^0}\Big[{1\over N}\sum_{k=1}^N\sup_{t_j \le s\le t_{j+1}} |\D X^{\btau, k}_s|\Big] &\le&   {1\over 3}\dbE^{\dbP^0}\Big[ \sup_{t_j\le s\le t_{j+1}}  \|\D \bX^\btau_s\|_1\Big] +C\dbE^{\dbP^0}\Big[\|\D X^{\btau}_{t_j}\|_1\Big]  \\
 &&+ {C\over N} \sum_{k=1}^N (1 + |x_k| +\|\bx\|_1 +|\D x_k|)|\D i_k|
\eeaa
Then one can easily get
\beaa
 \dbE^{\dbP^0}\Big[{1\over N}\sum_{k=1}^N\sup_{t \le s\le T} |\D X^{\btau, k}_s|\Big] \le   C\Big[\|\D \by\|_1+ (\|\bx\|_2+ \|\D \bx\|_2) \|\D \bi\|_2\Big],
\eeaa
which implies \reff{unif-est4} immediately.

\no{\bf Step 2.} In this step we prove the claimed regularity of $v_N$ in terms of $\by$. Fix $R>1$ and let $\tilde R>R$ be a large number which will be specified later.  Denote 
\beaa
D_{\tilde R} := [0, T]\times \{(x, m)\in \dbR^d\times \cP_2(\bS): |x| + \|m\|_2 \le \tilde R\}.
\eeaa
Note in particular that any bounded subset of $\cP_2(\bS)$ is tight in $\cP_1(\bS)$. Then $f$ and $g$ are uniformly continuous for $\cW_1$ on $D_{\tilde R}$, and let  $\tilde \rho_{\tilde R}$ denote their common modulus of continuity function. Then, by the linear growth of $f$ and $g$, 
\beaa
\label{tilderho}
\left.\ba{lll}
\dis |f(t,x,m)-f(t, \tilde x, \tilde m)| \le \tilde \rho_{\tilde R}\big(|\D x| + \cW_1(m, \tilde m)\big) \\
\dis \qq \qq \qq \qq + C(1+|x|+ \|m\|_1)\1_{\{|x|+\|m\|_2\ge \tilde R\}} + C(1+|\tilde x| + \|\tilde m\|_1) \1_{\{|\tilde x|+\|\tilde m\|_2\ge \tilde R\}}\Big]\\
\dis \qq \qq \qq \qq \qq \le \tilde \rho_{\tilde R}\big(|\D x| + \cW_1(m, \tilde m)\big) + {C\over \tilde R}(1+ |x|^2 + \|m\|_2^2 + |\tilde x|^2 + \|\tilde m\|_2^2);\\
\dis |g(m)-g(\tilde m)| \le \tilde \rho_{\tilde R}\big(\cW_1(m, \tilde m)\big)+ {C\over \tilde R}(1+  \|m\|_2^2 +  \|\tilde m\|_2^2).
\ea\right.
\eeaa 
We now fix $t\in [0, T]$, $\by, \tilde \by \in \bS^N$ with $\|\by\|_2, \|\tilde\by\|_2 \le R$.
By \reff{unif-estimate21} we have
\bea
\label{unif-est6}
&& \dis \dbE^{\dbP^0}\Big[\big|f(s, X^{\btau, k}_s, m^N(\bY^\btau_s)) - f(s, \tilde X^{\btau, k}_s, m^N(\tilde \bY^\btau_s))\big|\Big]\nonumber\\
&&\dis \qq \qq \qq \le \dbE^{\dbP^0}\Big[ \tilde \rho_{\tilde R}(|\D X^{\btau, k}_s| + \cW_1(m^N(\bY^\btau_s), m^N(\tilde \bY^\btau_s))\big)\nonumber \\
&&\dis \qq \qq \qq \qq + {C\over \tilde R} \big(1+ |X^{\btau, k}_s|^2 + \|m^N(\bY^\btau_s)\|_2^2 + |\tilde X^{\btau, k}_s|^2 + \|m^N(\tilde \bY^\btau_s)\|_2^2 \big)\Big]  \nonumber\\
&&\dis \qq \qq \qq \le \dbE^{\dbP^0}\Big[ \tilde \rho_{\tilde R}(|\D X^{\btau, k}_s| + \cW_1(m^N(\bY^\btau_s), m^N(\tilde \bY^\btau_s))\big) \Big] + {C(1+R^2)\over \tilde R};\\
 &&\dis \dbE^{\dbP^0}\Big[ \big|g(m^N(\bY^\btau_T)) - g(m^N(\tilde \bY^\btau_T))\big|\Big] \le  \dbE^{\dbP^0}\Big[ \tilde \rho_{\tilde R}( \cW_1(m^N(\bY^\btau_T), m^N(\tilde \bY^\btau_T))\big) \Big] + {C(1+R^2)\over \tilde R}.\nonumber
\eea
Then
\beaa
\dis \Big|J_N(t, \by, \btau)-J_N(t, \tilde \by, \btau)\Big| &\le& \dbE^{\dbP^0}\Big[{1\over N}\sum_{k=1}^N \int_t^{\t_k}  \tilde \rho_{\tilde R}(|\D X^{\btau, k}_s|+\cW_1(m^N(\bY^\btau_s), m^N(\tilde \bY^\btau_s))\big)ds\nonumber \\
&&\dis + \tilde \rho_{\tilde R}\big(\cW_1(m^N(\bY^\btau_T), m^N(\tilde \bY^\btau_T))\big) \Big] + {C(1+R^2)\over \tilde R}.
\eeaa
Now for any $\e>0$, choose $\tilde R = \tilde R_\e := {2C(1+R^2)\over \e}$.  By \reff{DItau} and \reff{unif-est4}, and observing that
\bea\label{ineq-modulus}
&&\dbE^{\dbP^0}\Big[\tilde \rho_{\tilde R_\e}\big(\cW_1(m^N(\bY^\btau_s), m^N(\tilde \bY^\btau_s))\big)\Big] \\ &&\le \tilde \rho_{\tilde R_\e}(\g) + \frac{1}{\sqrt{\g}}\sqrt{\dbE^{\dbP^0}\Big[\tilde \rho_{\tilde R_\e}^2\big(\cW_1(m^N(\bY^\btau_s), m^N(\tilde \bY^\btau_s))\big)\Big] \dbE^{\dbP^0}\Big[\cW_1(m^N(\bY^\btau_s), m^N(\tilde \bY^\btau_s))\big)\Big] \nonumber}
\eea
for all $\g > 0$ and $s \in [t,T]$, one can easily see that there exists $\d > 0$, depending on $\tilde \rho_{\tilde R_\e}$ and hence on $\e$, but not on $N$ or $\btau$, such that, whenever $\|\D\by\|_2 < \d$,
 \beaa
\dbE^{\dbP^0}\Big[{1\over N}\sum_{k=1}^N \! \int_t^{\t_k} \!\!\!\!\! \tilde \rho_{\tilde R_\e}(|\D X^{\btau, k}_s|+\cW_1(m^N(\bY^\btau_s), m^N(\tilde \bY^\btau_s))\big) ds  + \tilde \rho_{\tilde R_\e}\big(\cW_1(m^N(\bY^\btau_T), m^N(\tilde \bY^\btau_T))\big)\Big] \le {\e\over 2}.
\eeaa
Thus $\Big|J_N(t, \by, \btau)-J_N(t, \tilde \by, \btau)\Big| \le \e$.  By the arbitrariness of $\btau$, we obtain
\beaa
\Big|v_N(t, \by)-v_N(t, \tilde \by)\Big|\le \e,\q\mbox{whenever}\q \|\by\|_2, \|\tilde \by\|_2 \le R,~ \|\D \by\|_2 < \d.
\eeaa
This implies \reff{unif-estimate32} when $t=\tilde t$, that is, for a modulus of continuity function $\rho_R'$, 
\bea
\label{unif-est5}
\Big|v_N(t, \by)-v_N(t, \tilde \by)\Big|\le \rho_R'(\cW_2(m^N(\by), m^N(\tilde \by))),\q\mbox{whenever}\q \|\by\|_2, \|\tilde \by\|_2 \le R.
\eea

\no{\bf Step 3.} Fix $t<\tilde t$ and $\by \in \bS^N$ with $\|\by\|_2 \le R$.  For any $\btau\in \mathbfcal{T}_{\tilde t,T}^N$, by \reff{Nparticles}  we have
\beaa
&&\dis J_N(\tilde t, \by, \btau) - v_N(t, \by) \le J_N(\tilde t, \by, \btau) - J_N(t, \by, \btau) \nonumber\\
&&\dis =  J_N(\tilde t, \by, \btau)- \dbE^{\dbP^0}_{t, \by}\Big[ \frac 1 N \sum_{k=1}^N \int_t^{\tilde t} f(r,  X^{\btau,k}_r, m^N(\bY_r^\btau)) I_r^{\btau, k}dr + J_N(\tilde t, \bY^\btau_{\tilde t}, \btau) \Big]\nonumber\\
&&\dis\le J_N(\tilde t, \by, \btau)- \dbE^{\dbP^0}_{t, \by}\Big[ J_N(\tilde t, \bY^\btau_{\tilde t}, \btau) \Big] + \dbE^{\dbP^0}_{t, \by}\Big[ \frac C N \sum_{k=1}^N \int_t^{\tilde t} \big(1+ |X^{\btau,k}_r| + \|\bY_r^\btau\|_1\big)dr\Big]\nonumber\\
&&\dis \le J_N(\tilde t, \by, \btau)- \dbE^{\dbP^0}_{t, \by}\Big[ J_N(\tilde t, \bY^\btau_{\tilde t}, \btau) \Big] + C(1+R)\sqrt{\tilde t-t},
\eeaa
where the last inequality is due to \reff{unif-estimate11}. 
Again let  $\tilde R>R$ be a large number which will be specified later.  By \reff{unif-est5}, \reff{unif-estimate31}, and \reff{unif-estimate11} we have
\beaa
&&\dis J_N(\tilde t, \by, \btau) - \dbE^{\dbP^0}_{t, \by}\Big[ J_N(\tilde t, \bY^\btau_{\tilde t}, \btau) \Big]\\
&&\dis \qq \q \le \dbE^{\dbP^0}_{t, \by}\Big[ \rho'_{\tilde R}\big(\cW_2(m^N(\by), m^N(\bY^\btau_{\tilde t}))\big)  + \big(|J_N(\tilde t, \by, \btau)| +|J_N(\tilde t, \bY^\btau_{\tilde t}, \btau)|\big) \1_{\{\|\bY^\btau_{\tilde t}\|_2> \tilde R\}} \Big]\\
&&\dis \qq \q \le \dbE^{\dbP^0}_{t, \by}\Big[ \rho'_{\tilde R}\big(\|\bY^\btau_{\tilde t} - \by\|_2\big)  + {C\over \tilde R}\big(|J_N(\tilde t, \by, \btau)|^2 +|J_N(\tilde t, \bY^\btau_{\tilde t}, \btau)|^2 +\|\bY^\btau_{\tilde t}\|_2^2 \big) \Big]\\
&&\dis \qq \q \le \dbE^{\dbP^0}_{t, \by}\Big[ \rho'_{\tilde R}\big(\|\bY^\btau_{\tilde t} - \by\|_2\big)  + {C\over \tilde R}\big(\|\by\|_2^2 +\|\bY^\btau_{\tilde t}\|_2^2 \big) \Big]\\
&&\dis \qq \q \le \dbE^{\dbP^0}_{t, \by}\Big[ \rho'_{\tilde R}\big(\|\bY^\btau_{\tilde t} - \by\|_2\big)  \Big] + {C(1+R^2)\over \tilde R}.
\eeaa
Then
\beaa
 J_N(\tilde t, \by, \btau) - v_N(t, \by) \le \dbE^{\dbP^0}_{t, \by}\Big[ \rho'_{\tilde R}\big(\|\bY^\btau_{\tilde t} - \by\|_2\big)  \Big] + {C(1+R^2)\over \tilde R}+ C[1+R]\sqrt{\tilde t-t}.
\eeaa
Now for any $\e>0$, set $\tilde R = \tilde R_\e := {3C(1+R^2)\over \e}$, and $\d>0$ small enough. Then for $\tilde t-t\le \d$,
\beaa
&&\dis J_N(\tilde t, \by, \btau) - v_N(t, \by) \le \dbE^{\dbP^0}_{t, \by}\Big[ \rho'_{\tilde R}\big(\|\bY^\btau_{\tilde t} - \by\|_2\big)  \Big] +{2\e\over 3}.
\eeaa
Since $\btau\in \mathbfcal{T}_{\tilde t,T}^N$, we have $\t_k\ge \tilde t$ for all $k$, then $\bI^\btau_{\tilde t-} = \bi$. Thus, by \reff{unif-estimate12},
\beaa
\dbE^{\dbP^0}_{t, \by}\Big[\|\bY^\btau_{\tilde t} - \by\|_2\Big] = \dbE^{\dbP^0}_{t, \by}\Big[\|\bX^\btau_{\tilde t} - \bx\|_2\Big] \le C(1+R)(\tilde t-t)\le C(1+R)\d.
\eeaa
Using an estimate similar to \eqref{ineq-modulus}, one can then see that, for $\d>0$ sufficiently small, which may depend on $R$ and $\tilde R_\e$ and hence on $\e$, but not on $\btau$ or $N$, we have
\beaa
J_N(\tilde t, \by, \btau) - v_N(t, \by) \le \e,\q\mbox{whenever}\q \|\by\|_2 \le R,\q \tilde t-t\le \d.
\eeaa
Then, by the arbitrariness of $\btau\in \mathbfcal{T}_{\tilde t,T}^N$, we have
\bea
\label{unif-est7}
v_N(\tilde t, \by) - v_N(t, \by) \le \e,\q\mbox{whenever}\q \|\by\|_2 \le R,\q \tilde t-t\le \d.
\eea
On the other hand, for any $\btau\in \mathbfcal{T}_{t,T}^N$, note that $\btau \vee \tilde t\in \mathbfcal{T}_{\tilde t,T}^N$. Then similarly we have
\bea
\label{unif-est8}
 J_N(t, \by, \btau \vee \tilde t) - v_N(\tilde t, \by) \le   \e,\q\mbox{whenever}\q \|\by\|_2 \le R,\q \tilde t-t\le \d.
\eea
Note also that, 
\beaa
\cW_1\big(m^N(\bX_{\btau}^\btau), m^N(\bX_{\btau\vee \tilde t}^\btau)\big) \le \|\bX_{\btau}^\btau - \bX_{\btau\vee \tilde t}^\btau\|_1 \le C\sup_{t\le s\le \tilde t} \|\bX^\btau_s - \bx\|_1.
\eeaa
Then, by \reff{unif-estimate21}, \reff{unif-estimate22}, and the last estimate in \reff{unif-est6}, we have
\beaa
&&\dis  \Big|J_N(t, \by, \btau) - J_N(t, \by, \btau \vee \tilde t)\Big| \\
&&\dis = \Big|\dbE^{\dbP^0}_{t, \by}\Big[ \frac 1 N \sum_{k=1}^N \int_{\t_k}^{\t_k \vee \tilde t} f(r,  X^{\btau,k}_r, m^N(\bY_r^\btau)) dr + g\big(m^N(\bX_{\btau}^\btau) \big) - g\big(m^N(\bX_{\btau\vee \tilde t}^\btau) \big) \Big] \Big|\\
&&\dis \le  C\dbE^{\dbP^0}_{t, \by}\Big[ \frac 1 N \sum_{k=1}^N \int_t^{\tilde t} \big(1+ |X^{\btau,k}_r| + \|\bY_r^\btau\|_1\big) dr + \tilde \rho_{\tilde R}( \cW_1(m^N(\bX_{\btau}^\btau), m^N(\bX_{\btau\vee \tilde t}^\btau)\big) \Big] + {C(1+R^2)\over \tilde R}\\
&&\dis \le  C\dbE^{\dbP^0}_{t, \by}\Big[ \tilde \rho_{\tilde R}\big(C\sup_{t\le s\le \tilde t} \|\bX^\btau_s - \bx\|_1\big) \Big]  + {C(1+R^2)\over \tilde R} +
 C[1+R] \sqrt{\tilde t-t}.
\eeaa
This, together with \reff{unif-est8}, implies that, for a possibly smaller $\d$, 
\beaa
J_N(t, \by, \btau) - v_N(\tilde t, \by) \le  \e,\q\mbox{whenever}\q \|\by\|_2 \le R,\q \tilde t-t\le \d.
\eeaa
Since $\btau\in \mathbfcal{T}_{t,T}^N$ is arbitrary, we have
\beaa
v_N(t, \by) - v_N(\tilde t, \by) \le  \e,\q\mbox{whenever}\q \|\by\|_2 \le R,\q \tilde t-t\le \d.
\eeaa
Then by \reff{unif-est7} we obtain 
\beaa
\Big|v_N(t, \by) - v_N(\tilde t, \by)\Big| \le  \e,\q\mbox{whenever}\q \|\by\|_2 \le R,\q \tilde t-t\le \d.
\eeaa
This, together with \reff{unif-est5} and \reff{Birkhoff}, proves \reff{unif-estimate32}.
\qed

\subsection{Construction of an optimal stopping strategy}
In this subsection we come back to the general case in Subsection \ref{sect-general}. We first recall  \reff{Ai}, and note that, by \reff{vNmon},
\bea
\label{vNmon2}
\max_{\bi'<\bi} v^N(t, \bx, \bi') = \max_{k\in A(\bi)} v^N(t, \bx, \bi^{-k})
\eea
We shall use the above relation to construct an optimal   $\hat\btau \in \mathbfcal{T}_{0,T}^N$ for the problem \eqref{Nparticles}. We also refer to \cite[Section 4.3]{TTZ} for some heuristic discussions on the optimal stopping policy on the mean field problem \reff{weakoptstop}.

Fix $\by \in \bS^N$ and assume for simplicity $t=0$. We set
\bea
\label{btau}
\hat \btau := \btau^{(|A(\bi)|)},
\eea
where $\btau^{(j)}\in \mathbfcal{T}_{0,T}^N$, $j=0,\cds, |A(\bi)|$, are constructed as follows.

{\bf Step 0 (Initialization).} Set  $\tau^{(0)}_k := T\1_{A(\bi)}(k) $.  That is, $\bY^{\btau^{(0)}} = \bar \bY$ is the unstopped process. Moreover, set $A_0 := A(\bi)$ and $\t^*_0:=0$. 

We then stop the particles in $A(\bi)$ recursively: for $j=0,\cds, |A(\bi)|-1$:

{\bf Step 1.} The next particle is stopped at: denoting $\bI^A := (1_A(1), \cds, 1_A(N))$ for $A\subset [N]$,
\bea
\label{t*next}
\t^*_{j+1} := \inf\Big\{s \ge \t^*_{j} : v^N(s, \bY_s^{\btau^{(j)}}) = \underset{\bi' < \bI^{A_j}}{\max} \ v^N(s, \bX_s^{\btau^{(j)}}, \bi')\Big\}.
\eea
We remark that it is possible that $\t^*_{j+1} = \t^*_j$, namely two (or more) particles are stopped at the same time.

{\bf Step 2.} By \reff{vNmon2}, we may pick an index for the particle stopped at $\t^*_{j+1}$:  
\bea
\label{k*next}
\k^*_{j+1}:= \min\Big\{k\in A_j:   v^N(s, \bX_s^{\btau^{(j)}}, (\bI^{A_j})^{-k}) = v^N(s, \bY_s^{\btau^{(j)}})\Big\}, \q A_{j+1}:= A_j \backslash \{\k^*_{j+1}\},
\eea
where $\bI^{-k}$ is defined in \reff{bi-l}, and $A_{j+1} \in \cF^N_{\t^*_{j+1}}$ is random. Note that there might be multiple $k$ satisfying the above requirement and serving for our purpose. We choose the smallest one just for convenience. 

{\bf Step 3.} We then update the vector of the stopping times: 
\bea
\label{btaunext}
 \t^{(j+1)}_{\k^*_{j+1}} := \t^*_{j+1}, \q\mbox{and}\q  \tau^{(j+1)}_k := \tau^{(j)}_k,~\mbox{for all}~ k\neq \k^*_{j+1}. \nonumber
 \eea

\begin{thm}\label{verifN}
Under Assumption \ref{assum-N}, 
the $\hat\btau$ constructed in \reff{btau} is an optimal stopping policy for the problem \eqref{Nparticles}.
\end{thm}
\proof We prove by induction on $|A(\bi)|$. When $|A(\bi)|=0$, all particles are already stopped. Then it is trivial that $\hat\btau = \btau^{(0)} = (0,\cds, 0)$ is optimal. Assume the result is true for all $\bi'$ with $|A(\bi')| \le n$, and now consider $\bi$ such that $|A(\bi)| = n+1$. 

Recall \reff{Ntostd} and note that $\bY^{\btau^{(0)}}=\bar \bY$, $\bI^{A_0} = \bi$. Then it follows from the standard optimal stopping theory that $\t^*_1$ constructed in \reff{t*next} is an optimal stopping strategy for \eqref{Ntostd} with $t=0$. That is,
\beaa
v^N(0,\by) &=& \dbE_{0,\by}\Big[\frac 1 N \sum_{k=1}^N \int_0^{\t^*_1} f^k(r, \bar \bX_r)\bar I_r^kdr + \max_{\bi' < \bi} v^N(\t^*_1, \bar \bX_{\t^*_1}, \bi') \Big]\\
&=&\dbE_{0,\by}\Big[\frac 1 N \sum_{k\in A_0} \int_0^{\t^*_1} f^k(r, \bX_r^{\btau^{(0)}}) dr + v^N(\t^*_1, \bX^{\btau^{(0)}}_{\t^*_1}, \bI^{A_1}) \Big],
\eeaa
where the last equality follows from \reff{k*next}. Our construction also induces that
\beaa
\bX^{\hat \t}_r =  \bX_r^{\btau^{(j)}},\q \t^*_j \le r\le \t^*_{j+1},\q\mbox{and}\q \bI^{\hat \btau}_r =  \bI^{A_j},\q \t^*_j \le r< \t^*_{j+1}.
\eeaa
Then
\beaa
v^N(0,\by) =\dbE_{0,\by}\Big[\frac 1 N \sum_{k=1}^N \int_0^{\t^*_1} f^k(r, \bX_r^{\hat \btau}) I^{\hat \btau, k}_r dr + v^N(\t^*_1, \bX^{\hat \btau}_{\t^*_1}, \bI^{A_1}) \Big].
\eeaa
Since $|A_1|= n$, by induction assumption we have
\beaa
v^N(\t^*_1, \bX^{\hat \btau}_{\t^*_1}, \bI^{A_1}) = \dbE_{\t^*_1, \bX^{\hat \btau}_{\t^*_1}, \bI^{A_1}}\Big[\frac 1 N \sum_{k=1}^N \int_{\t^*_1}^T f^k(r, \bX_r^{\hat \btau}) I^{\hat \btau, k}_r dr + g_N(\bX^{\hat \btau}_T) \Big],
\eeaa
thus proving the optimality of $\hat\btau$:
\beaa
v^N(0,\by) =\dbE_{0,\by}\Big[\frac 1 N \sum_{k=1}^N \int_0^T f^k(r, \bX_r^{\hat \btau}) I^{\hat \btau, k}_r dr + g_N(\bX^{\hat \btau}_T) \Big].
\eeaa

\vspace{-1cm}
\qed

\section{Appendix}
\label{sect-appendix}
In this appendix, we present a few technical proofs. 

\ms
\no{\bf Proof of Lemma \ref{pure-mixed}.} Let $m' \preceq m$ with corresponding $p(x)$ as in \reff{order}. For any $n\ge 1$, let $\{O^n_k\}_{k\ge 1}$ be a partition of $\dbR^d$ such that, for each $k\ge 1$,  $|x-\tilde x|\le {1\over n}$ for all $x, \tilde x\in O^n_k$. By the required continuity of $m(dx, 1)$, one may easily construct $A^n_k \subset O^n_k$ such that
\bea
\label{Ank}
m(A^n_k, 1) = \int_{O^n_k}p(x) m(dx, 1). \nonumber
\eea
Set $A^n:= \cup_{k\ge 1} A^n_k$ and fix an $x^n_k\in A^n_k$. Then, for any $\f: \bS\to \dbR$ which is Lipschitz continuous with Lipschitz constant $1$, we have: denoting $\D \f(x) := \f(x,1) - \f(x,0)$, 
\beaa
&&\dis \int_\bS \f(\by) m^{A^n}(dy) - \int_\bS \f(\by) m'(dy) \\
&&\dis = \int_{A^n} \f(x, 1) m(dx,1) +  \int_{(A^n)^c} \f(x, 0) m(dx,1) + \int_{\dbR^d} \f(x, 0) m(dx,0) \\
&&\dis \q - \int_{\dbR^d} \f(x,1) p(x) m(dx,1) - \int_{\dbR^d}\Big(\f(x,0) [1-p(x)] m(dx,1) + \f(x,0) m(dx,0)\Big)\\
&&\dis = \int_{A^n} \D\f(x)m(dx,1)  - \int_{\dbR^d} \D\f(x) p(x) m(dx,1)\\
&&\dis = \sum_{k\ge 1} \Big( \int_{A^n_k} \D\f(x) m(dx,1)  - \int_{O^n_k} \D\f(x) p(x) m(dx,1)\Big)\\
&&\dis = \sum_{k\ge 1} \Big( \int_{A^n_k} (\D\f(x)-\D \f(x^n_k)) m(dx,1) - \int_{O^n_k} (\D\f(x)-\D \f(x^n_k)) p(x) m(dx,1)\Big)\\
&&\dis \le {C\over n}\sum_{k\ge 1} \Big( \int_{A^n_k} m(dx,1)  + \int_{O^n_k} p(x) m(dx,1)\Big) \le {C\over n}.
\eeaa
This implies that $\cW_1(m^{A_n}, m') \le {C\over n}$, and thus $m^{A_n}\to m'$ weakly as $n\to\infty$. Moreover, it is obvious that $\{m^{A^n}\}_{n\ge 1}$ are uniformly square integrable, then it follows from \cite[Theorem 5.5]{CarDel} that $\lim_{n\to\infty} \cW_2(m^{A_n}, m') =0$.
\qed

\ms
\no{\bf Proof of Lemma \ref{lem-derivative}.} 
First, let $h\in \dbR^d$ and denote $\by +_k h\in \bS^N$ as the perturbation of $\by$ by replacing $x_k$ with $x_k + h$. Set $m_h^{l, k} := l m^N(\by +_k h) + (1-l) m^N(\by)$. Since $\f \in C_2^{1,2}(\bQ_0)$ and $i_k = 1$, we have
\beaa
\dis\phi(s, \by +_k h) - \phi(s, \by)&=&  \f(s, m^N(\by +_k h)) - \f(s, m^N(\by)) 
\\
&=&\dis  \int_0^1 \int_\bS \d_m \f\big(s, m_h^{l, k}, y\big)d(m^N(\by +_k h)-m^N(\by))(y)dl 
\\
&=&\dis \frac 1 N \int_0^1 \Big[\d_m \f\big(s, m_h^{l, k}, (x_k + h, 1)\big) - \d_m \f\big(s, m_h^{l, k}, (x_k,1)\big)\Big]dl
\\
&=& \frac 1 N \int_0^1 \int_0^1 h\cd \pa_x \d_m \f_1\big(s, m_h^{l, k}, x_k + \tilde l h\big)d\tilde l dl.
\eeaa
Note that $x_k +  \tilde l h \to x_k$ and  $m_h^{l, k} \underset{h \to 0}{\longrightarrow} m^N(\by)$, uniformly in $\tilde l, l\in [0,1]$, as $|h|\to 0$. Then
\beaa
\phi(s, \by +_k h) - \phi(s, \by) = h \cd {1\over N} \pa_x \d_m \f_1\big(s, m^N(\by), x_k\big) + o(|h|).
\eeaa
This implies immediately  the first equality in \reff{connection-derivatives}.

Similarly, one can show that
$$
\dis \pa_x \d_m \big\{\f_1(s, m^N(\by+_k h), x_k) - \f_1(s, m^N(\by), x_k) \big\}
=
{1\over N} \pa_{\m \m}^2 \f_{1,1}(s, m^N(\by), x_k, x_k) h + o(|h|).
$$
Then
\begin{align*}
\pa_{x_k}\big\{\phi(s, \by +_k h) -\phi(s, \by)\big\}
&=\frac 1 N \pa_x \d_m \big\{\f_1(s, m^N(\by+_k h), x_k+h) - \f_1(s, m^N(\by), x_k)\big\}
\\
&\hspace{-15mm}=
\frac 1 N \pa_{xx}^2 \d_m \f_1(s, m^N(\by), x_k) h + \frac{1}{N^2} \pa_{\m \m}^2 \f_{1,1}(s, m^N(\by), x_k, x_k) h + o(|h|).
\end{align*}
This implies the second equality in \reff{connection-derivatives}.
\qed

{\small
\bibliographystyle{plain}
\bibliography{references2}}

\end{document}